\documentclass[a4paper,11pt,reqno]{amsart}

\usepackage[headings]{fullpage}

\usepackage{amssymb}
\usepackage{dsfont}

\usepackage[utf8]{inputenc} 
\usepackage[T1]{fontenc} 

\usepackage{hyperref}

\usepackage[nobysame,alphabetic,initials,msc-links]{amsrefs}

\DefineSimpleKey{bib}{how}

\renewcommand{\eprint}[1]{#1}
\BibSpec{misc}{%
  +{}{\PrintAuthors}  {author}
  +{,}{ \textit}      {title}
  +{,}{ }             {how}
  +{}{ \parenthesize} {date}
  +{,} { available at \eprint}        {eprint}
  +{,}{ available at \url}{url}
  +{,}{ }             {note}
  +{.}{}              {transition}
}

\usepackage{wasysym}
\newcommand{\circt}%
{\mathbin{%
\mathchoice
{\ooalign{$\ocircle$\cr\hidewidth\raise-.15ex\hbox{$\scriptstyle\top\mkern2.05mu$}\cr}}
{\ooalign{$\ocircle$\cr\hidewidth\raise-.15ex\hbox{$\scriptstyle\top\mkern2.05mu$}\cr}}
{\ooalign{$\scriptstyle\ocircle$\cr\hidewidth\raise-.12ex\hbox{$\scriptscriptstyle\top\mkern1mu$}\cr}}
{\ooalign{$\scriptstyle\ocircle$\cr\hidewidth\raise-.12ex\hbox{$\scriptscriptstyle\top\mkern1mu$}\cr}}
}}

\DeclareFontFamily{U}{matha}{\hyphenchar\font45}
\DeclareFontShape{U}{matha}{m}{n}{
      <5> <6> <7> <8> <9> <10> gen * matha
      <10.95> matha10 <12> <14.4> <17.28> <20.74> <24.88> matha12
      }{}
\DeclareSymbolFont{matha}{U}{matha}{m}{n}
\DeclareFontFamily{U}{mathb}{\hyphenchar\font45}
\DeclareFontShape{U}{mathb}{m}{n}{
      <5> <6> <7> <8> <9> <10> gen * mathb
      <10.95> mathb10 <12> <14.4> <17.28> <20.74> <24.88> mathb12
      }{}
\DeclareSymbolFont{mathb}{U}{mathb}{m}{n}
\DeclareMathSymbol{\ovoid}{3}{matha}{"6C}
\DeclareMathSymbol{\boxvoid}{2}{mathb}{"6C}

\mathchardef\mhyph="2D

\numberwithin{equation}{section}

\newtheorem{theorem}{Theorem}[section]
\newtheorem{corollary}[theorem]{Corollary}
\newtheorem{lemma}[theorem]{Lemma}
\newtheorem{proposition}[theorem]{Proposition}
\theoremstyle{remark}
\newtheorem{conjecture}[theorem]{Conjecture}
\newtheorem{remark}[theorem]{Remark}

\theoremstyle{definition}
\newtheorem{definition}[theorem]{Definition}

\newcommand\bp{\begin{proof}}
\newcommand\ep{\end{proof}}

\DeclareMathOperator{\Ad}{Ad}
\DeclareMathOperator{\ad}{ad}
\DeclareMathOperator{\Rep}{Rep}

\DeclareMathOperator\Dhat{\hat\Delta}
\DeclareMathOperator\wt{\operatorname{wt}}

\newcommand{\C}{{\mathbb C}}

\newcommand{\Z}{{\mathbb Z}}

\newcommand\T{{\mathbb T}}

\newcommand{\A}{{\mathcal A}}
\newcommand{\B}{{\mathcal B}}

\newcommand\CC{{\mathcal C}}

\newcommand{\F}{{\mathcal F}}
\newcommand\HH{\mathcal H}

\newcommand{\K}{{\mathcal K}}

\newcommand\OO{\mathcal O}
\newcommand\RR{\mathcal R}

\newcommand\TT{\mathcal T}
\newcommand\U{\mathcal U}

\newcommand\g{\mathfrak g}
\newcommand\h{\mathfrak h}

\newcommand\eps{\varepsilon}

\begin{document}

\title{Cartan subproduct systems}

\author{Suvrajit Bhattacharjee}
\address{Wydział Matematyki, Uniwersytet w Białymstoku, K. Ciołkowskiego 1M, 15-245 Białystok, Poland}
\email{s.bhattacharjee@uwb.edu.pl}

\author{Olof Giselsson}
\address{Department of Mathematical Sciences, Chalmers University of Technology and
the University of Gothenburg, Gothenburg SE-412 96, Sweden}
\email{olofgi@chalmers.se}

\author{Sergey Neshveyev}
\address{Department of Mathematics, University of Oslo, PB 1053 Blindern, 0316 Oslo, Norway}
\email{sergeyn@math.uio.no}

\thanks{The research  was supported by the NFR project 300837 ``Quantum Symmetry''. S.B. was also supported by the National Science Centre, Poland, through the WEAVE-UNISONO grant no. 2023/05/Y/ST1/00046.}

\begin{abstract}
Given a semisimple compact Lie group $G$ and a nonzero dominant integral weight~$\lambda$, the highest weight $G_q$-modules $V_{n\lambda}$ form a subproduct system of finite dimensional Hilbert spaces. Using a conjectural asymptotic behavior of Clebsch--Gordan coefficients we identify the corresponding Cuntz--Pimsner algebras with quantized algebras of functions on homogeneous spaces of $G$. We also show that the gauge-invariant part of the Toeplitz algebra provides a model for convergence of full matrix algebras to quantum flag manifolds, complementing and generalizing results of Landsman and Rieffel for $q=1$ and results of Vaes--Vergnioux in the rank one case for $q\ne1$.

We verify our conjecture on Clebsch--Gordan coefficients for $G=SU(n)$ and all weights that are either regular or multiples of the fundamental weight $\omega_1$. For $\lambda=\omega_1$, we also provide a detailed description of the Toeplitz and Cuntz--Pimsner algebras, generalizing results of Arveson on symmetric subproduct systems.
\end{abstract}

\date{December 15, 2025; revised August 23, 2026}

\maketitle

\section*{Introduction}

Operator algebras associated to product systems of C$^*$-correspondences form a large class of algebras that includes, in particular, Cuntz algebras, graph algebras and semigroup C$^*$-algebras. An even larger class of algebras arises from subproduct systems~\cites{MR2646788,MR2608451}. The algebras are defined in similar ways in both settings: they are generated by creation operators $S_\xi$ on Fock-type spaces. What, however, makes the general case of subproduct systems more complicated is that there are no obvious relations involving both creation and annihilation operators, apart from the basic inequality $\sum_iS_{\xi_i}S_{\xi_i}^*\le1$ when $\sum_i\theta_{\xi_i,\xi_i}\le1$. At present there are no general methods of obtaining relations in these C$^*$-algebras, describing their ideal structure or computing their K-theory.

Recently the third author together with Erik Habbestad introduced a class of so called Temperley--Lieb subproduct systems~\cites{MR4705666,MR4919591} having a large quantum symmetry, in the sense that the associated Cuntz--Pimsner algebras carry a natural ergodic action of a compact quantum group. The idea of exploiting symmetries of the construction was of course not new. Already in the foundational work of Arveson~\cite{arveson} the action of the unitary group was used in a crucial way to describe the corresponding Cuntz--Pimsner algebras, and later a study of equivariant subproduct systems was explicitly initiated by Andersson~\cite{Anders}, followed by the work of Arici and Kaad~\cite{MR4451552}. The novelty of~\cite{MR4705666} lied in bringing modern tools of quantum group theory, such as monoidal equivalence and a categorical formulation of the Baum--Connes conjecture, to analyze the associated algebras, describe their generators and relations and compute the K-theory.
Temperley--Lieb subproduct systems are defined by particular quadratic polynomials in noncommuting variables. More general quadratic polynomials have been now analyzed by Aiello, Del Vecchio and Rossi~\cite{MR4985585}, and recently a systematic study of quadratic subproduct systems has been undertaken by Arici and Ge~\cite{AG}.

In the present paper we consider another natural class of equivariant subproduct systems that includes Arveson's symmetric subproduct systems and, up to monoidal equivalence, the Temperley--Lieb subproduct systems. Namely, take a simply connected semisimple compact Lie group $G$. Consider its $q$-deformation $G_q$ and take an irreducible highest weight $G_q$-module~$V_\lambda$. Then, for all $m,n\ge0$, the module $V_{m\lambda}\otimes V_{n\lambda}$ contains a unique submodule isomorphic to~$V_{(m+n)\lambda}$, which is called the Cartan component of the tensor product. This way we get a subproduct system $(V_{n\lambda})_{n\ge0}$ and the associated Cuntz--Pimsner algebra $\OO_{\lambda,q}$.

An attempt to analyze the gauge-invariant subalgebra of $\OO_{\lambda,q}$ has been made in~\cite{Anders}, but unfortunately~\cite{Anders} contains a number of gaps that seem difficult to fill (the most critical of which is the proof of the overly optimistic \cite{Anders}*{Theorem 4.9}). At the heart of our analysis of $\OO_{\lambda,q}$ lies a derivation of approximate mixed commutation relations between the creation and annihilation operators. In fact, finding a method to do this when there are neither simple Jones--Wenzl type formulas for the projections $V_\lambda^{\otimes n}\to V_{n\lambda}$ as in~\cites{MR4705666,MR4919591}, nor hopes to find explicit expressions for the annihilation operators by a brute force computation as in~\cite{arveson}, was one of the main motivations for our work.

The method we found relies on a conjectural property of the Clebsch--Gordan coefficients that says roughly that when $q\ge1$ and $n$ is large, then every highest weight vector in $V_\lambda\otimes V_{n\lambda}$ is close to a vector of the form $\xi\otimes \xi_{n\lambda}$, where $\xi_{n\lambda}$ is a highest weight vector in $V_{n\lambda}$, and every vector of this form with $\xi\perp\xi_\lambda$ is close to be orthogonal to $V_{(n+1)\lambda}$, see Section~\ref{ssec:conj} for the precise formulation. For $q=1$ the first property follows from the main result of~\cite{MR1600423}, but in our opinion the proof in~\cite{MR1600423} is incomplete and does not address key analytical difficulties. We haven't been able to prove the conjecture in full generality, but we verified it for $G=SU(N)$ and all weights $\lambda$ that are either regular or multiples of the fundamental weight $\omega_1$ (corresponding to the standard representation of $SU(N)$ on $\C^N$).

\smallskip

The paper is organized as follows. In the preliminary Section~\ref{sec:cqg} we collect definitions and constructions related to subproduct systems and quantum groups. The main new result here is a sufficient condition for ergodicity of the action of a compact quantum group on the Cuntz--Pimsner algebra of an equivariant subproduct system.

In Section~\ref{sec:q-sym}, as a warm up for the general case, we study the Toeplitz and Cuntz--Pimsner algebras of $q$-symmetric subproduct systems, which correspond to the case $G=SU(N)$ and $\lambda=\omega_1$. In this case everything can be computed explicitly. The compact quantum group~$SU_q(N)$ plays an important role in the analysis, but we use very little beyond its definition. Our results here generalize some of the results of Arveson for $q=1$~\cite{arveson}.

Section~\ref{sec:cartan} contains our main results. The strategy of our analysis of $\OO_{\lambda,q}$ is as follows (it is presented in Section~\ref{sec:cartan} in the opposite order). We use the conjectural asymptotic behavior of the Clebsch--Gordan coefficients to obtain enough relations in $\OO_{\lambda,q}$ to be able to conclude that the action of $G_q$ on $\OO_{\lambda,q}$ is ergodic. We then construct a character on $\OO_{\lambda,q}$. This allows us to define an injective $G_q$-equivariant homomorphism $\OO_{\lambda,q}\to C(G_q)$. It is then a relatively routine task to identify its image.

In Section~\ref{sec:gauge} we obtain additional results on the gauge-invariant subalgebra $\TT^{(0)}_{\lambda,q}$ of the Toeplitz algebra. We first show that $\TT^{(0)}_{\lambda,q}$ is the section algebra of a continuous field of C$^*$-algebras over $\Z_+\cup\{\infty\}$ with fibers $\B(V_{n\lambda})$ at $n\in\Z_+$ and $C(K^S_q\backslash G_q)$ at $\infty$, where $K^S_q\backslash G_q$ is a quantized flag manifold. This can be interpreted as a convergence $\B(V_{n\lambda})\to C(K^S_q\backslash G_q)$. Results of this type have a long history. For $q=1$ the existence of such a continuous field structure was first proved by Landsman~\cite{MR1656992} using properties of the Berezin quantization. Later Rieffel showed ~\cite{MR2055928} that the convergence $\B(V_{n\lambda})\to C(K^S\backslash G)$ can be understood rigorously within the theory of quantum metric spaces, see also~\cite{MR2195335} for an alternative approach and~\cites{MR2200270,MR4426737} for some related work in the quantum setting. Complementing these results, we show that the continuous field structure defined by $\TT^{(0)}_\lambda$ coincides with the one introduced in~\cite{MR1656992}.

For $q\ne1$ another interpretation of convergence $\B(V_{n\lambda})\to C(K^S_q\backslash G_q)$ can be given in terms of compactifications of discrete quantum spaces. For $G_q=SU_q(2)$ (and, more generally, for the monoidally equivalent free orthogonal quantum groups) such compactifications were constructed by Vaes and Vergnioux~\cite{MR2355067} using properties of the Jones--Wenzl projections. Generalizing results of~\cite{MR4705666}, we show that the conjectural asymptotic behaviour of the Clebsch--Gordan coefficients is enough to construct similar compactifications in the higher rank case. Moreover, our analysis of $\TT^{(0)}_{\lambda,q}$ allows us to achieve two goals in one strike by proving that the compactifications exist and the corresponding boundaries are the quantized flag manifolds~$K^S_q\backslash G_q$.

Let us finish by noting that for a number of results of the paper it would be more natural to consider subproduct systems over the monoid $P_+$ of dominant integral weights rather than over $\Z_+$. Although we briefly take this point of view in Section~\ref{ssec:star-comm}, a detailed analysis of such systems will be carried out elsewhere.

\bigskip

\section{Compact quantum groups and equivariant subproduct systems}\label{sec:cqg}

\subsection{Compact quantum groups}
Let us fix the notation and recall some basic notions, see~\cite{neshveyev-tuset-book} for more details. By a compact quantum group $G$ we mean a Hopf $*$-algebra $(\C[G],\Delta)$ that is generated by matrix coefficients of finite dimensional unitary right comodules. We have a one-to-one correspondence between such comodules and the finite dimensional unitary representations of $G$, that is, unitaries $U\in \B(H_U)\otimes\C[G]$ such that $(\iota\otimes\Delta)(U)=U_{12}U_{13}$. Namely, the right comodule structure on $H_U$ is given by
$$
\delta_U\colon H_U\to H_U\otimes\C[G],\quad \delta_U(\xi)=U(\xi\otimes1).
$$
The tensor product $U\otimes V$ of two unitary representations (denoted also by $U\circt V$ or $U\times V$) is defined by $U_{13}V_{23}$.

Denote by $h$ the Haar state on $\C[G]$ and by $L^2(G)$ the corresponding GNS-space. We view $\C[G]$ as a subalgebra of $\B(L^2(G))$. The (reduced) C$^*$-algebra $C(G)$ of continuous functions on~$G$ is defined as the norm closure of $\C[G]$.

Consider the $*$-algebra $\U(G)=\C[G]^*$ dual to the coalgebra $(\C[G],\Delta)$, with the involution $\omega^*(x)=\overline{\omega(S(x)^*)}$. Every finite dimensional unitary representation $U$ of $G$ defines a $*$-representation
$$
\pi_U\colon\U(G)\to \B(H_U),\quad \pi_U(\omega)=(\iota\otimes\omega)(U).
$$
We usually suppress $\pi_U$ and write $\omega\xi$ instead of $\pi_U(\omega)\xi$.

\smallskip

Assume now that $A$ is a C$^*$-algebra and $\alpha\colon A\to A\otimes C(G)$ is a $*$-homomorphism such that $(\alpha\otimes\iota)\alpha=(\iota\otimes\Delta)\alpha$.
We say that $\alpha$ is a (right) action of $G$ on $A$, or that $A$ is a $G$-C$^*$-algebra, if either of the following equivalent conditions is satisfied:
\begin{enumerate}
  \item[(1)] the linear space $\alpha(A)(1\otimes C(G))$ is a dense subspace of $A\otimes C(G)$ (the Podle\'s condition);
  \item[(2)] there is a dense $*$-subalgebra $\A\subset A$ on which $\alpha$ defines a coaction of the Hopf algebra $(\C[G],\Delta)$, that is, $\alpha(\A)\subset\A\otimes_{\mathrm{alg}}\C[G]$ and $(\iota\otimes\eps)\alpha(a)=a$ for all $a\in\A$,
\end{enumerate}
where $\otimes_{\mathrm{alg}}$ denotes the purely algebraic tensor product.
Then the largest subalgebra as in (2) is given by
$$
\A=\operatorname{span}\{(\iota\otimes h)(\alpha(a)(1\otimes x)): a\in A,\ x\in\C[G]\},
$$
its elements are called regular.

An action is called reduced if $\alpha$ is injective, or equivalently, by faithfulness of $h$ on $C(G)$, if the conditional expectation
$$
E=(\iota\otimes h)\alpha\colon A\to A^G
$$
is faithful. Here $A^G$ denotes the fixed point subalgebra $\{a\in A:\alpha(a)=a\otimes1\}$. If $G$ is coamenable, that is, the counit $\eps$ on $\C[G]\subset C(G)$ is bounded, then $(\iota\otimes\eps)\alpha=\iota$ on~$A$, so all actions of $G$ are reduced. An action on a unital C$^*$-algebra $A$ is called ergodic if $A^G=\C1$.

\subsection{Subproduct systems}
Recall, following~\citelist{\cite{MR2608451}\cite{MR2646788}}, that a subproduct system  $\HH$ of finite dimensional Hilbert spaces (over the additive monoid $\Z_+$) is a sequence of Hilbert spaces $(H_n)^\infty_{n=0}$ together with isometries $w_{k,l}\colon H_{k+l}\to H_k\otimes H_l$ such that
$$
\dim H_0=1,\qquad \dim H_1=m<\infty,\qquad (w_{k,l}\otimes1)w_{k+l,n}=(1\otimes w_{l,n})w_{k,l+n}.
$$
By~\cite{MR2608451}*{Lemma~6.1}, we can assume that $H_0=\C$, $H_{k+l}\subset H_k\otimes H_l$ and the isometries $w_{k,l}$ are simply the embedding maps. The subproduct systems satisfying these stronger properties are called standard. For such subproduct systems we denote by $f_n$ the projection $H_1^{\otimes n}\to H_n$.

Given a subproduct system $\HH=(H_n)^\infty_{n=0}$, the associated Fock space is defined by
$$
\F_\HH=\bigoplus^\infty_{n=0}H_n.
$$
For every $\xi\in H_1$, we define a creation operator
$$
S_\xi\colon\F_\HH\to\F_\HH\quad\text{by}\quad S_\xi\zeta=w_{1,n}^*(\xi\otimes\zeta)\quad\text{for}\quad \zeta\in H_n.
$$
We will often fix an orthonormal basis $(\xi_i)^N_{i=1}$ in $H_1$ and write $S_i$ for $S_{\xi_i}$. The Toeplitz algebra~$\TT_\HH$ of $\HH$ is defined as the unital C$^*$-algebra generated by $S_1,\dots,S_N$.

If $\HH$ is standard and $H=H_1$, it is convenient to identify $\F_\HH$ with a subspace of the full Fock space
$
\F(H)=\bigoplus^\infty_{n=0}H^{\otimes n}.
$
Consider the operators $T_i\colon\F(H)\to\F(H)$, $T_i\zeta=\xi_i\otimes\zeta$, and the projection
$
e_\HH\colon \F(H)\to\F_\HH.
$
Then $S_i=e_\HH T_i|_{\F_\HH}$ and $S_i^*=T_i^*|_{\F_\HH}$.
Since $1-\sum^N_{i=1}T_iT_i^*$ is the projection onto $H^{\otimes 0}=\C$, we then get
\begin{equation}\label{eq:e-0}
e_0=1-\sum^N_{i=1}S_iS_i^*,
\end{equation}
where $e_0$ is the projection $\F_\HH\to H_0$. From this one concludes that $\K(\F_\HH)\subset\TT_\HH$. The Cuntz--Pimsner algebra~\cite{MR2949219} of $\HH$ is defined by
$$
\OO_\HH=\TT_\HH/\K(\F_\HH).
$$
We will denote the images of $S_i$ and $S_\xi$ in $\OO_\HH$ by $s_i$ and $s_\xi$, resp.

\smallskip

Once we fix an orthonormal basis $(\xi_i)^N_{i=1}$ in $H$, it is often convenient to identify the tensor algebra~$T(H)$ with the algebra $\C\langle X_1,\dots,X_N\rangle$ of polynomials in $N$ noncommuting variables. 
By~\cite{MR2608451}*{Proposition~7.2}, there is a one-to-one correspondence between the standard subproduct systems with $H_1=H$ and the homogeneous ideals $I$ in $\C\langle X_1,\dots,X_N\rangle$ such that the degree one homogeneous component $I_1$ of $I$ is zero. Namely, given such an ideal $I$, we define
$$
H_n=I_n^\perp\subset H^{\otimes n}.
$$
Then $\F_\HH$ can be identified with a completion of the algebra $\C\langle X_1,\dots,X_N\rangle/I$. In this picture the creation operators $S_i$ are the unique bounded extensions of the operators of multiplication on the left by $X_i$ on this algebra, and it becomes obvious that $p(S_1,\dots,S_N)=0$ for all $p\in I$.

In general it seems difficult to say anything about relations between the operators $S_i$ and $S_j^*$, apart from the inequality $\sum_iS_iS_i^*\le1$ that follows from~\eqref{eq:e-0}. It is not even clear when the following properties hold.

\begin{definition}
Given a subproduct system $\HH$, we say that $\TT_\HH$ admits \emph{normal ordering} if the elements $S_{i_1}\dots S_{i_n}S^*_{j_1}\dots S^*_{j_m}$ span a dense subspace of $\TT_\HH$. Similarly, we say that $\OO_\HH$ admits normal ordering if the elements $s_{i_1}\dots s_{i_n}s^*_{j_1}\dots s^*_{j_m}$ span a dense subspace of $\OO_\HH$.
\end{definition}

We remark that the widely open Arveson--Douglas conjecture (see, e.g.,~\cite{Hartz}*{Section~3.5}) predicts that $\OO_\HH$ is commutative if the ideal $I$ corresponding to $\HH$ contains the commutators $[X_i,X_j]$. When $\OO_\HH$ is commutative, it obviously admits normal ordering.

\subsection{Equivariant subproduct systems}
Assume $G$ is a compact quantum group and $\HH=(H_n)^\infty_{n=0}$ is a subproduct system of finite dimensional unitary $G$-modules. By this we mean that we are given unitary representations $U_n$ of $G$ on $H_n$ and the structure maps $H_{k+l}\to H_k\otimes H_l$ are $G$-intertwiners. In this case we say also that $(H_n)^\infty_{n=0}$ is a $G$-equivariant subproduct system.

Given such a system, we have a unitary representation
$$
U_\HH=\bigoplus^\infty_{n=0} U_n\in M(\K(\F_\HH)\otimes C(G))
$$
of $G$ on $\F_\HH$. It defines reduced right actions of $G$ on $\K(\F_\HH)$ and $\TT_\HH$ by
$$
\alpha(T)=U_\HH(T\otimes1)U_\HH^*.
$$
The action on $\TT_\HH$ passes to a possibly nonreduced action on $\OO_\HH$. We have the following sufficient condition for ergodicity of this action.

\begin{proposition}\label{prop:normal}
Assume $G$ is a compact quantum group and $\HH=(H_n)_{n\ge0}$ is a $G$-equivariant subproduct system of finite dimensional Hilbert spaces. Assume that
\begin{enumerate}
    \item the $G$-modules $H_n$ are simple and pairwise nonisomorphic;
    \item the Cuntz--Pimsner algebra $\OO_\HH$ admits normal ordering.
\end{enumerate}
Then the action of $G$ on $\OO_\HH$ is ergodic.
\end{proposition}

We remind that ergodicity means that the fixed point algebra $\OO_\HH^G$ is trivial.

\bp
Consider the linear maps $S^{(n)}\colon H_1^{\otimes n}\to\OO_\HH$ defined by $S^{(n)}(\xi_1\otimes\dots\otimes\xi_n)=s_{\xi_1}\dots s_{\xi_n}$. They are $G$-equivariant and have the property $S^{(n)}=S^{(n)}f_n$, where $f_n$ is the projection $H_1^{\otimes n}\to H_n$. More generally, consider the linear maps $S^{(n,m)}\colon H_1^{\otimes n}\otimes \bar H_1^{\otimes m}\to\OO_\HH$ defined by
$$
S^{(n,m)}(\xi_1\otimes\dots\otimes\xi_n\otimes\bar\zeta_1\otimes\dots\otimes\bar\zeta_m)=s_{\xi_1}\dots s_{\xi_n}s^*_{\zeta_1}\dots s^*_{\zeta_m}.
$$
Then $S^{(n,m)}=S^{(n,m)}(f_n\otimes\bar f_m)$, where $\bar f_m\colon \bar H_1^{\otimes m}\to\bar H_m$ is the coisometry defined by
$$
\bar f_m(\bar\zeta_1\otimes\dots\otimes\bar\zeta_m)=\overline{f_m(\zeta_m\otimes\dots\otimes\zeta_1)}
$$
and we identify $\bar H_m$ with a subspace of $\bar H_1^{\otimes m}$ using the adjoint map, so that $\bar f_m$ becomes a projection. The maps $S^{(n,m)}$ become equivariant if we consider the contragredient action of~$G$ on $\bar H_1$, that is, the corresponding $\U(G)$-module structure is given by
$\omega\bar\xi=\overline{\hat S(\omega)^*\xi}$, where $\hat S(\omega)=\omega S$ and $S$ in the antipode on $(\C[G],\Delta)$.

Consider the conditional expectation $E\colon\OO_\HH\to\OO_\HH^G$ defined by $E(a)=(\iota\otimes h)\alpha(a)$, where $\alpha\colon \OO_\HH\to\OO_\HH\otimes C(G)$ is the action of $G$ and $h$ is the Haar state on $C(G)$. By the normal ordering assumption (ii) we then conclude that the elements of the form $S^{(n,m)}(v)$, where $v\in H_n\otimes\bar H_m$ is a $G$-invariant vector, span a dense subspace of $\OO_\HH^G$. By assumption (i) there are no such nonzero vectors $v$ for $n\ne m$, while for $n=m$ the space of such vectors is one-dimensional. Therefore, for every $n\ge1$, the space $E(S^{(n,n)}(H_1^{\otimes n}\otimes \bar H_1^{\otimes n}))=E(S^{(n,n)}(H_n\otimes \bar H_n))$ is at most one-dimensional. This space does contain the unit, since once we fix an orthonormal basis in~$H_1$, we get $\sum_{i_1,\dots,i_n}s_{i_1}\dots s_{i_n}s_{i_n}^*\dots s_{i_1}^*=1$ by~\eqref{eq:e-0}. Hence $\OO_\HH^G=\C1$.
\ep

\subsection{Drinfeld--Jimbo deformations of semisimple compact Lie groups}\label{ssec:DJ}
Let $G$ be a simply connected semisimple compact Lie group, $\g$ its
complexified Lie algebra. The universal enveloping algebra $U\g$ is a Hopf $*$-algebra with involution such that the real Lie algebra of~$G$ consists of skew-adjoint elements. Fix a nondegenerate symmetric $\ad$-invariant form on $\g$ such that its restriction to the real Lie algebra of $G$ is negative definite. Let $\h\subset\g$ be the Cartan subalgebra defined by a maximal torus $T$ in $G$.
When $G$ is simple, we always normalize the symmetric form so that for the dual form on $\h^*$ and every short root $\alpha\in\Delta$ we have $(\alpha,\alpha)=2$.

For every root $\alpha$ put $d_\alpha=(\alpha,\alpha)/2$. Let $H_\alpha\in\h$ be the element corresponding to the coroot $\alpha^\vee=2\alpha/(\alpha,\alpha)$ under the identification $\h\cong\h^*$.
Fix a system $\Pi=\{\alpha_1,\dots,\alpha_r\}$ of simple roots. For every positive root $\alpha\in\Delta_+$ choose $E_\alpha\in\g_\alpha$ such that $(E_\alpha,E_\alpha^*)=d_\alpha^{-1}$, and put $F_\alpha=E_\alpha^*\in \g_{-\alpha}$, so that $[E_\alpha,F_\alpha]=H_\alpha$. We write $E_i,F_i,H_i$ for $E_{\alpha_i},F_{\alpha_i}, H_{\alpha_i}$, resp.  Denote by $\omega_1,\dots,\omega_r$ the fundamental weights, so $\omega_i(H_j)=\delta_{ij}$.

\smallskip

Fix a number $q>0$. If $q=1$, we put $U_1\g=U\g$. For $q\ne1$ the quantized universal
enveloping algebra~$U_q\g$ is generated by elements $E_i$, $F_i$, $K_i$,
$K_i^{-1}$, $1\le i\le r$, satisfying the relations
$$
K_iK_i^{-1}=K_i^{-1}K_i=1,\ \ K_iK_j=K_jK_i,\ \
K_iE_jK_i^{-1}=q_i^{a_{ij}}E_j,\ \
K_iF_jK_i^{-1}=q_i^{-a_{ij}}F_j,
$$
$$
E_iF_j-F_jE_i=\delta_{ij}\frac{K_i-K_i^{-1}}{q_i-q_i^{-1}},
$$
$$
\sum^{1-a_{ij}}_{k=0}(-1)^k\begin{bmatrix}1-a_{ij}\\
k\end{bmatrix}_{q_i} E^k_iE_jE^{1-a_{ij}-k}_i=0,\ \
\sum^{1-a_{ij}}_{k=0}(-1)^k\begin{bmatrix}1-a_{ij}\\
k\end{bmatrix}_{q_i} F^k_iF_jF^{1-a_{ij}-k}_i=0,
$$
where $a_{ij}=(\alpha_i^\vee,\alpha_j)=\alpha_j(H_i)$, $\displaystyle\begin{bmatrix}m\\
k\end{bmatrix}_{q_i}=\frac{[m]_{q_i}!}{[k]_{q_i}![m-k]_{q_i}!}$,
$[m]_{q_i}!=[m]_{q_i}[m-1]_{q_i}\dots [1]_{q_i}$,
$\displaystyle[n]_{q_i}=\frac{q_i^n-q_i^{-n}}{q_i-q_i^{-1}}$,
$q_i=q^{d_i}$ and $d_i=d_{\alpha_i}$. This is a Hopf $*$-algebra with coproduct $\Dhat_q$ and
counit $\hat\eps_q$ defined by
$$
\Dhat_q(K_i)=K_i\otimes K_i,\ \
\Dhat_q(E_i)=E_i\otimes1+ K_i\otimes E_i,\ \
\Dhat_q(F_i)=F_i\otimes K_i^{-1}+1\otimes F_i,
$$
$$
\hat\eps_q(E_i)=\hat\eps_q(F_i)=0,\ \ \hat\eps_q(K_i)=1,
$$
and with involution given by $K_i^*=K_i$, $E_i^*=F_iK_i$, $F_i^*=K_i^{-1}E_i$.

If $V$ is a finite dimensional $U_q\g$-module and $\lambda\in
P$ is an integral weight, denote by $V(\lambda)$ the
space of vectors $v\in V$ of weight $\lambda$, so that
$K_iv=q^{(\lambda,\alpha_i)}v=q_i^{\lambda(i)}v$ for all $i$, where $\lambda(i)=(\lambda,\alpha_i^\vee)=\lambda(H_i)$. Recall that $V$ is called admissible if
$V=\oplus_{\lambda\in P}V(\lambda)$. Every such module is unitarizable, so we can always assume that $V$ is a Hilbert space and the representation $U_q\g\to \B(V)$ is a $*$-homomorphism. We denote by $\CC_q(\g)$ the tensor category of finite dimensional admissible unitary $U_q\g$-modules. This category is semisimple and its simple objects are highest weight modules $V_\lambda$ for dominant integral weights $\lambda\in P_+$.

Denote by $\C[G_q]\subset (U_q\g)^*$ the Hopf $*$-algebra of matrix coefficients of finite dimensional admissible $U_q\g$-modules. This defines a compact quantum group with representation category~$\CC_q(\g)$. The fusion rules in this category and weight decompositions of highest weight modules are the same for all $q>0$.

Recall that $\U(G_q)$ denotes the dual $*$-algebra of $\C[G_q]$. The maximal torus $T\subset G$ can be considered as a subgroup of the unitary group of $\U(G_q)$ and its complexification $T_\C$ can be considered as a subgroup of the group of invertible elements of $\U(G_q)$: if $X\in\mathfrak h$, then for any admissible $U_q\g$-module $V$ and $\lambda\in P$ the element $\exp( X)\in T_\C$ acts on~$V(\lambda)$ as multiplication by $e^{\lambda(X)}$. Under this embedding $T_\C\hookrightarrow\U(G_q)$, we have $K_i=q_i^{H_i}$ for $i=1,\dots,r$.

We will need the following classical result, see, e.g., \cite{Ku}*{Proposition 3.2 and Corollary~3.4}.

\begin{theorem}\label{thm:mult}
For any dominant integral weights $\lambda,\mu,\nu\in P_+$, the multiplicity $m^\nu_{\lambda,\mu}$ of $V_\nu$ in $V_\lambda\otimes V_\mu$ is not larger than $\dim V_\lambda(\nu-\mu)$. The equality holds if $\lambda'(i)+\mu(i)\ge-1$ for all $i$ and all weights $\lambda'$ in the weight decomposition of~$V_\lambda$.
\end{theorem}

The category $\Rep G_q$ is braided, with braiding $\sigma=\Sigma\RR$, where $\Sigma\colon V\otimes W\to W\otimes V$ is the flip map and $\RR\in \U(G_q\times G_q)$ is the universal $R$-matrix uniquely determined by the property that if $\zeta\in V$ is a lowest weight vector of weight $\lambda$ and $\xi\in W$ is a highest weight vector of weight $\mu$ for some finite dimensional simple $G_q$-modules $V$ and $W$, then
\begin{equation} \label{eq:ermat0}
\RR(\zeta\otimes\xi)=q^{(\lambda,\mu)}\zeta\otimes\xi.
\end{equation}
Explicitly, the $R$-matrix has the form
\begin{equation} \label{eq:ermat1}
\RR=q^{\sum_{i,j}(B^{-1})_{ij}H_i\otimes H_j}\prod_{\alpha\in\Delta_+}\exp_{q_\alpha}((1-q_\alpha^{-2})F_\alpha\otimes E_\alpha),
\end{equation}
where $B$ is the matrix $((\alpha_i^\vee,\alpha_j^\vee))_{i,j}$, the product is taken with respect to a particular order on the set~$\Delta_+$ of positive roots,
$$
\exp_q(\omega)=\sum^\infty_{n=0}q^{n(n+1)/2}\frac{\omega^n}{[n]_q!},
$$
and $E_\alpha$ (resp. $F_\alpha$) is a polynomial in $K_i^{\pm1}$ and $E_i$ (resp. $F_i$) such that $K_iE_\alpha=q^{(\alpha_i,\alpha)}E_\alpha K_i$ (resp. $K_iF_\alpha=q^{-(\alpha_i,\alpha)}F_\alpha K_i$), so that $E_\alpha V(\lambda)\subset V(\lambda+\alpha)$ (resp., $F_\alpha V(\lambda)\subset V(\lambda-\alpha)$), see \cite{CP}*{Theorem~8.3.9} for details (in the conventions of \cite{CP} we have $q=e^h$, $K_i=e^{-d_ihH_i}$, $E_i=X_i^-$ and $F_i=X_i^+$).

From~\eqref{eq:ermat1} we see that~\eqref{eq:ermat0} holds for more general pairs of weight vectors. Namely, writing $\wt(\xi)$ for the weight of $\xi$, we have
\begin{equation} \label{eq:ermat2}
\RR(\zeta\otimes\xi)=q^{(\wt(\zeta),\wt(\xi))}\zeta\otimes\xi
\end{equation}
for all weight vectors $\zeta$ and $\xi$ such that either $\zeta$ is a lowest weight vector or $\xi$ is a highest weight vector.

\smallskip

Next let us recall a description of the Poisson--Lie subgroups of $G$ with respect to the standard $r$-matrix and their $q$-deformations~\cite{MR2914062}. Take a subset $S$ of $\Pi$. Denote by $\tilde K^S$ the closed connected subgroup of $G$ such that its complexified Lie algebra $\tilde\g_S$ is generated by the elements $E_i$ and $F_i$ with $\alpha_i\in S$, so
$$
\tilde\g_S=\operatorname{span}\{H_i\mid \alpha_i\in S\}\oplus\bigoplus_{\alpha\in\Delta_S}\g_\alpha,
$$
where $\Delta_S$ is the set of roots that lie in the group generated by $\alpha_i\in S$. Denote by $P(S^c)$ the subgroup of the weight lattice $P$ generated by the fundamental weights $\omega_i$ with $\alpha_i\in S^c=\Pi\setminus S$. Let $L$ be a subgroup of $P(S^c)$. Identifying $P$ with the dual group of the maximal torus $T$, denote by $T_L$ the annihilator of $L$ in $T$. Since $T$ normalizes $\tilde K^S$, the group $K^{S,L}$ generated by $\tilde K^S$ and $T_L$ is a closed subgroup of $G$, and its complexified Lie algebra is
$$
\g_{S,L}=\h_L\oplus\bigoplus_{\alpha\in\Delta_S}\g_\alpha,
$$
where $\h_L\subset\h$ is the annihilator of $L\subset\h^*$. Note that if $L=P(S^c)$ then $K^{S,L}$ is the group $\tilde K^S$. If $L=0$, we write $K^S$ for $K^{S,L}$. Then $K^S=G\cap P_S$, where $P_S\subset G_\C$ is the parabolic subgroup corresponding to $S$, and $\tilde K^S$ is the semisimple part of $K^S$.

For $q>0$, $q\ne1$, let $\U(K^{S,L}_q)$ be the $\sigma(\U(G_q),\C[G_q])$-closed subalgebra of $\U(G_q)$ generated by $T_L$ and $E_i$, $F_i$ with $\alpha_i\in S$. In other words, an element $\omega\in\U(G_q)$ belongs to $\U(K^{S,L}_q)$ if and only if for every finite dimensional admissible $U_q\g$-module $V$ the operator of the action by~$\omega$ on~$V$ lies in the algebra generated by $T_L$ and $E_i$, $F_i$ with $\alpha_i\in S$. Denote by $\C[K^{S,L}_q]\subset \U(K^{S,L}_q)^*$ the Hopf $*$-algebra that is the image of $\C[G_q]$ under the restriction map $\U(G_q)^*\to\U(K^{S,L}_q)^*$.

We will need the following refinement of the first part of Theorem~\ref{thm:mult}.

\begin{proposition}\label{prop:mult}
Assume we are given  $\lambda,\mu,\nu\in P_+$. Let $S=\{\alpha\in\Pi:(\nu,\alpha)=(\mu,\alpha)=0\}$ and $L=\Z(\nu-\mu)$. Then the multiplicity of $V_\nu$ in $V_\lambda\otimes V_\mu$ is not larger than the dimension of the subspace of vectors in $V_\lambda(\nu-\mu)$ fixed by $K^{S,L}_q$, that is, vectors $\xi$ such that $\omega\xi=\hat\eps_q(\omega)\xi$ for all $\omega\in\U(K^{S,L}_q)$.
\end{proposition}

\bp
Consider the contragredient module $\bar V_\mu$, so that $\omega\bar\eta=\overline{\hat S_q(\omega)^*\eta}$ for $\omega\in U_q\g$ and $\eta\in V_\mu$, where $\hat S_q$ is the antipode on $(U_q\g,\Dhat_q)$. Let $\zeta\in\bar V_\mu$ be a lowest weight vector. It has weight $-\mu$. Let $\xi\in V_\nu$ be a highest weight vector, of weight $\nu$. Then the vector $\xi\otimes\zeta\in V_\nu\otimes\bar V_\mu$ is cyclic. Since both $\xi$ and $\zeta$ are killed by $E_i$ and $F_i$ for $\alpha_i\in S$, the same is true for $\xi\otimes\zeta$. Also, this vector has weight $\nu-\mu$. Hence the multiplicity of $V_\nu$ in $V_\lambda\otimes V_\mu$, which is equal to the dimension of the space of morphisms $V_\nu\otimes\bar V_\mu\to V_\lambda$, is not larger than the dimension of the subspace of vectors in $V_\lambda(\nu-\mu)$ fixed by $K^{S,L}_q$.
\ep

By construction we have an epimorphism $\pi\colon\C[G_q]\to\C[K^{S,L}_q]$ of Hopf $*$-algebras. Put
$$
\C[G_q/K^{S,L}_q]=\{a\in \C[G_q]\mid (\iota\otimes\pi)\Delta_q(a)=a\otimes 1\},
$$
where $\Delta_q$ is the comultiplication on $\C[G_q]$. Equivalently, $a\in \C[G_q/K^{S,L}_q]$ if and only if $(\iota\otimes \omega)\Delta_q(a)=\hat\eps_q(\omega)a$ for all $\omega\in\U(K^{S,L}_q)$. Denote by $C(G_q/K^{S,L}_q)$ the norm-closure of  $\C[G_q/K^{S,L}_q]$ in $C(G_q)$. The C$^*$-algebra $C(G_q/K^{S,L}_q)$ can also be defined as a universal C$^*$-completion of~$\C[G_q/K^{S,L}_q]$, see~\cite{MR2914062}*{Remark~2.1}. In a similar way we can define $C(K^{S,L}_q\backslash G_q)$.

For $\lambda\in P_+$ and $\xi,\zeta\in V_\lambda$ denote by $C^\lambda_{\zeta,\xi}\in \C[G_q]$ the matrix coefficient $(\cdot \,\xi,\zeta)$. Then $\C[G_q/K^{S,L}_q]$ is the linear span of elements $C^\lambda_{\zeta,\xi}$ such that $\lambda\in P_+$, $\zeta\in V_\lambda$ and $\xi\in V_\lambda$ is fixed by~$K^{S,L}_q$.

\bigskip

\section{A case study: \texorpdfstring{$q$}{q}-symmetric subproduct systems}\label{sec:q-sym}

Fix a natural number $N\ge 2$ and a real number $q>0$. In this section we study the subproduct system defined by the homogeneous ideal of $\C\langle X_1,\dots,X_N\rangle$ with generators $X_iX_j-qX_jX_i$ for $i<j$.

\subsection{A \texorpdfstring{$q$}{q}-analogue of the Drury--Arveson space}

The corresponding quotient of the algebra $\C\langle X_1,\dots,X_N\rangle$ is the quantized algebra of polynomials $\C_q[e_1\dots,e_N]$, where the generators $e_i$ satisfy the relations
$$
e_ie_j=qe_je_i, \quad i < j.
$$
The Fock space $H^2_{N,q}$ obtained by completing $\C_q[e_1,\dots,e_N]$ is a $q$-analogue of the Drury--Arveson space $H^2_N$. The space $H^2_{N,q}$ was introduced in a greater generality in~\cite{BB}. The definition of the scalar product in~\cite{BB} was however ad hoc, so it is not immediately obvious that it agrees with the one defined using subproduct systems. Our first goal is to verify that this is nevertheless true.

In fact, we will show a bit more. In the case $q=1$ it is known that the scalar product on $\C[e_1,\dots,e_N]$ is completely determined by the following properties, see~\cite{Hartz}*{Section~2.2}: (i) homogeneous polynomials of different degrees are orthogonal to each other; (ii) the scalar product is invariant under the action of $SU(N)$; (iii) $\|e_1^k\|=1$ for all $k\ge0$. We will see that the same is true for all $q>0$ if $SU(N)$ is replaced by the compact quantum group $SU_q(N)$.

It will be more convenient to work with the quantized universal enveloping algebra $U_q\mathfrak{sl}_N$. Recall that the standard representation of $U_q\mathfrak{sl}_N$ on $\C^N=\C e_1+\dots+\C e_N$ is given by
$$
E_ie_j=\delta_{i+1,j}q^{1/2}e_{j-1},\qquad
F_ie_j=\delta_{i,j}q^{-1/2}e_{j+1},\qquad K_ie_j=q^{\delta_{i,j}-\delta_{i+1,j}}e_j
$$
for $i=1,\dots,N-1$ and $j=1,\dots,N$. 
Using the formulas for the coproduct on $U_q\mathfrak{sl}_N$ it is not difficult to check that the subspace of $\C^N\otimes\C^N$ spanned by the vectors $e_i\otimes e_j-qe_j\otimes e_i$ for $i<j$ is $U_q\mathfrak{sl}_N$-invariant. It follows that we get a representation of the algebra $U_q\mathfrak{sl}_N$ on $\C_q[e_1,\dots,e_N]$. This representation is $*$-preserving with respect to the scalar product on $\C_q[e_1,\dots,e_N]\subset H^2_{N,q}$, since the representation on $\C^N$ is $*$-preserving. This is an equivalent form of $SU_q(N)$-invariance of the scalar product.

\begin{proposition}[{cf.~\cite{arveson}*{Lemma~3.8}, \cite{BB}*{Section~1}}]\label{prop:inner-products}
For all $N\ge2$ and $q>0$, the monomials $e_1^{d_1}\dots e_N^{d_N}$ ($d_1,\dots,d_N\in\Z_+$) form an orthogonal basis in $H^2_{N,q}$, and
\begin{equation}\label{eq:inner-products}
\lVert e_{1}^{d_{1}}\cdots e_{N}^{d_{N}} \rVert^{2}=q^{D}\frac{[d_{1}]_{q}!\cdots [d_{N}]_{q}!}{[d_{1}+\cdots+d_{N}]_{q}!},
\end{equation}
where $D=\sum_{1 \leq i < j \leq N}d_{i}d_{j}$.
\end{proposition}

\bp
It is obvious from the commutation relations that the monomials $e_1^{d_1}\dots e_N^{d_N}$ span the entire space $\C_q[e_1,\dots,e_N]$. By the definition of the Fock space $H^2_{N,q}$ it is also clear that the monomials of different degrees are mutually orthogonal. Different monomials $e_1^{d_1}\dots e_N^{d_N}$ of the same degree are also orthogonal, since they have different weights with respect to the action of~$U_q\mathfrak{sl}_N$. Therefore we only need to prove~\eqref{eq:inner-products}.

Observe that in $\C\langle X_1,\dots,X_N\rangle$ the monomial $X_1^k$ is orthogonal to any noncommutative polynomial of the form $P(X_1,\dots,X_N)(X_iX_j-q X_jX_i)Q(X_1,\dots,X_N)$, $i<j$. It follows that $\|e_1^k\|=1$. We will deduce~\eqref{eq:inner-products} from this using the action of $U_q\mathfrak{sl}_N$.

Let us first consider the case $N=2$. An easy induction argument yields the following formulas for the action of $U_{q}\mathfrak{sl}_{2}$ on the monomials in $\C_q[e_1,e_2]$:
\[
E(e_{1}^{m}e_{2}^{n})=q^{m-n+3/2}[n]_{q}e_{1}^{m+1}e_{2}^{n-1},\qquad
F(e_{1}^{m}e_{2}^{n})=q^{-m+n+1/2}[m]_{q}e_{1}^{m-1}e_{2}^{n+1}
\]
and $K(e_1^me_2^n)=q^{m-n}e_1^me_2^n$. As $F^*=K^{-1}E$, we then get
$$
( e_{1}^{m}e_{2}^{n},e_{1}^{m}e_{2}^{n} )
=\frac{q^{m-n+3/2}}{[m+1]_q}(Fe_{1}^{m+1}e_{2}^{n-1},e_{1}^{m}e_{2}^{n})
=q^{m-n+1}\frac{[n]_{q}}{[m+1]_{q}}(e_{1}^{m+1}e_{2}^{n-1},e_{1}^{m+1}e_{2}^{n-1}).
$$
From this we get by induction on $n$ that
\begin{equation*}\label{eq:inner-products-su-2}
\| e_{1}^{m}e_{2}^{n} \|^2=q^{mn}\frac{[m]_{q}![n]_{q}!}{[m+n]_{q}!}\|e_1^{m+n}\|^2
=q^{mn}\frac{[m]_{q}![n]_{q}!}{[m+n]_{q}!}.
\end{equation*}
This proves~\eqref{eq:inner-products} for $N=2$.

Consider now the general case. Observe that the copy of $U_q\mathfrak{sl}_2$ in $U_q\mathfrak{sl}_N$ generated by $E_{N-1}$, $F_{N-1}$ and $K_{N-1}$ acts trivially on $\C_q[e_1,\dots,e_{N-2}]$. The same argument as in the case $N=2$ shows then that for any $p\in\C_q[e_1,\dots,e_{N-2}]$ we have
$$
\lVert p e_{N-1}^me_N^n\rVert^{2}
=q^{mn}\frac{[m]_{q}![n]_{q}!}{[m+n]_{q}!}\lVert p e_{N-1}^{m+n}\rVert^{2}.
$$
Then we use the copy of $U_q\mathfrak{sl}_2$ in $U_q\mathfrak{sl}_N$ generated by $E_{N-2}$, $F_{N-2}$ and $K_{N-2}$, and so on. A simple induction on the length of a monomial gives~\eqref{eq:inner-products}.
\ep

\subsection{The Cuntz--Pimsner algebra}

Let us write $S_i$ for the creation operators $S_{e_i}$ on $H^2_{N,q}$. Thus,
\[
S_{i}(e_{1}^{d_{1}}\cdots e_{N}^{d_{N}})=e_{i}e_{1}^{d_{1}}\cdots e_{N}^{d_{N}}=q^{-D_{i-1}}e_{1}^{d_{1}}\cdots e_{i}^{d_{i}+1}\cdots e_{N}^{d_{N}},
\]
where $D_{i-1}=\sum_{j < i}d_{j}$. Let us also put $D_N=d_1+\dots+d_N$. Using Proposition~\ref{prop:inner-products} we get the following formula for the adjoint operators.

\begin{lemma}
For all $i=1,\dots,N$ and $d_1,\dots,d_N\in \Z_+$, we have
$$
S_{i}^{*}(e_{1}^{d_{1}}\cdots e_{N}^{d_{N}})=q^{D_N-D_{i}}\frac{[d_{i}]_{q}}{[D_N]_{q}}e_{1}^{d_{1}}\cdots e_{i}^{d_{i}-1}\cdots e_{N}^{d_{N}}.
$$
\end{lemma}

From this we will get the following result.

\begin{lemma}\label{lem:q-arveson-relations}
On the subspace $H_n\subset \C_q[e_1,\dots,e_N]$ spanned by the monomials of degree $n$, we have:
$$
S_i^*S_j=\frac{[n]_q}{[n+1]_q}S_jS_i^*\quad (i\ne j),
$$
$$
S_{i}^{*}S_{i}=q\frac{[n]_{q}}{[n+1]_{q}}S_{i}S_{i}^{*}+(q-q^{-1})\frac{[n]_{q}}{[n+1]_{q}}\sum^N_{j=i+1}S_{j}S_{j}^{*}+\frac{q^{-n}}{[n+1]_{q}}\mathrm{id}.
$$
\end{lemma}

\bp
The first relation is a simple consequence of the previous lemma, we will concentrate on the second one. The relation is obvious for $n=0$, so let us assume that $n\ge1$. Given $d_1,\dots,d_N\in\Z_+$ with $d_1+\dots+d_N=n$, we have
\begin{align*}
S_{i}^{*}S_{i}(e_{1}^{d_{1}}\cdots e_{N}^{d_{N}})&=q^{n-D_{i}-D_{i-1}}\frac{[d_{i}+1]_{q}}{[n+1]_{q}}e_{1}^{d_{1}}\cdots e_{N}^{d_{N}},\\
S_{j}S_{j}^{*}(e_{1}^{d_{1}}\cdots e_{N}^{d_{N}})&=q^{n-D_{j}-D_{j-1}}\frac{[d_{j}]_{q}}{[n]_{q}}e_{1}^{d_{1}}\cdots e_{N}^{d_{N}}.
\end{align*}
Therefore all we need to prove is the following identity:
\begin{equation*}
q^{n-D_{i}-D_{i-1}}[d_{i}+1]_{q}=q^{n+1-D_{i}-D_{i-1}}[d_i]_{q}+(q-q^{-1})\sum^N_{j=i+1}q^{n-D_{j}-D_{j-1}}[d_j]_{q}+q^{-n}.
\end{equation*}
Using that $[d_i+1]_q=q[d_i]_q+q^{-d_i}$ we can rewrite this as
\begin{equation*}
q^{n-2D_i}=(q-q^{-1})\sum^N_{j=i+1}q^{n-D_{j}-D_{j-1}}[d_j]_{q}+q^{-n}.
\end{equation*}
This, in turn, is easy to verify by a downward induction on $i$.
\ep

Denote by $\TT_{N,q}$ the C$^*$-subalgebra of $\B(H^2_{N,q})$ generated by the operators $S_i$ ($1\le i\le N$) and by $\OO_{N,q}$ the corresponding Cuntz--Pimsner algebra $\TT_{N,q}/\K(H^2_{N,q})$. Recall that we write $s_i$ for the image of $S_i$ in $\OO_{N,q}$. The following theorem generalizes \cite{arveson}*{Theorem~5.7}.

\begin{theorem}\label{thm:arveson}
For every $N\ge2$ and $q>0$, the Cuntz--Pimsner algebra $\OO_{N,q}$ is a universal unital C$^*$-algebra with generators $s_1,\dots,s_N$ satisfying the following relations:
\begin{enumerate}
    \item[(1)] if $q\ge1$, then
$$
\sum^N_{i=1}s_is_i^*=1,\qquad s_is_j=qs_js_i\quad (i<j),\qquad s_i^*s_j=q^{-1}s_js_i^*\quad (i\ne j),
$$
$$
s_{i}^{*}s_{i}=s_{i}s_{i}^{*}+(1-q^{-2})\sum^N_{j=i+1}s_{j}s_{j}^{*};
$$
\item[(2)] if $q\le 1$, then
$$
\sum^N_{i=1}s_is_i^*=1,\qquad s_is_j=qs_js_i\quad (i<j),\qquad s_i^*s_j=qs_js_i^*\quad (i\ne j),
$$
$$
s_{i}^{*}s_{i}=s_{i}s_{i}^{*}+(1-q^2)\sum^{i-1}_{j=1}s_{j}s_{j}^{*}.
$$
\end{enumerate}
\end{theorem}

\bp
We will only prove (1). Part (2) is proved similarly; it can also be deduced from (1) using the isomorphism $\C_q[e_1,\dots,e_N]\cong \C_{q^{-1}}[e_1,\dots,e_N]$, $e_i\mapsto e_{N-i+1}$.

The first identity in (1) holds for all subproduct systems, the second is an immediate consequence of the relations in $\C_q[e_1,\dots,e_N]$, the remaining two follow from the previous lemma. These relations describe the algebra of continuous functions on the quantum homogeneous space $S^{2N-1}_q=SU_q(N-1)\backslash SU_q(N)$, see~\cite{MR1086447}, where $SU(N-1)$ is embedded into $SU(N)$ as a right lower corner. Thus, we have a surjective $SU_q(N)$-equivariant homomorphism
$C(S^{2N-1}_q)\to\OO_{N,q}$. It must be an isomorphism, since $C(S^{2N-1}_q)$ has no nontrivial $SU_q(N)$-equivariant quotients.
\ep

Note that the relations in (2) describe the C$^*$-algebra $C(SU_q(N-1)\backslash SU_q(N))$, where $SU(N-1)$ is embedded into $SU(N)$ as a left upper corner. Both algebras $C(SU_q(N-1)\backslash SU_q(N))$ are well-defined for all $q>0$ and are isomorphic as C$^*$-algebras, which is obvious for $q=1$, while for $q\ne1$ this follows from the fact that all these algebras are isomorphic to one graph C$^*$-algebra obtained by formally letting $q^{-1}=0$ in (1) (for this one has to write the relation $s_is_j=qs_js_i$ as $q^{-1}s_is_j=s_js_i$) or $q=0$ in (2), see~\cite{MR1942860}. At the same time, for each $q\ne1$, they are not isomorphic as $SU_q(N)$-C$^*$-algebras.

\begin{remark}
Using Lemma~\ref{lem:q-arveson-relations} and Theorem~\ref{thm:arveson} it is possible to describe the Toeplitz algebras~$\TT_{N,q}$ in terms of generators and relations similarly to \cite{MR4705666}*{Theorem~3.4 and Remark~3.5}.
\end{remark}

\bigskip

\section{Cartan subproduct systems}\label{sec:cartan}

Throughout this section we fix a simply connected semisimple compact Lie group $G$ and $q>0$. For any dominant integral weights $\lambda,\mu\in P_+$, the $G_q$-module $V_\lambda\otimes V_\mu$ contains a unique copy of $V_{\lambda+\mu}$, which is called the \emph{Cartan component} of the tensor product. It follows that for every $\lambda\in P_+$ we get a subproduct system $(V_{n\lambda})_{n\ge0}$, which we call a \emph{Cartan subproduct system}; a~more formal definition will be given in Section~\ref{ssec:comm}. Our goal is to describe the corresponding Cuntz--Pimsner algebra $\OO_{\lambda,q}$.

\subsection{Stabilizers of highest weight vectors}

For $q>0$ and every $\lambda\in P_+$, fix a unit highest weight vector $\xi_\lambda\in V_\lambda$ and a lowest weight unit vector $\zeta_\lambda\in V_\lambda$. Note that $\wt(\zeta_\lambda)=w_0\lambda$, where~$w_0$ is the longest element of the Weyl group. Define an involution on $P$ by letting $\bar\lambda=-w_0\lambda$. We have the following result.

\begin{proposition}
For $\lambda\in P_+$, define
$$
S=\{\alpha\in\Pi: (\lambda,\alpha)=0\},\quad L=\Z\lambda.
$$
Then the following properties hold:
\begin{enumerate}
  \item[(1)] for $q=1$, the stabilizer $G^\lambda$ of $\xi_\lambda\in V_\lambda$ in $G$ coincides with $K^{S,L}$;
  \item[(2)] for each $q>0$, the $*$-algebra $\C[G_q/G_q^\lambda]$ is generated by the elements $C^\lambda_{\eta,\xi_\lambda}$ ($\eta\in V_\lambda$);
  \item[(3)] for each $q>0$, the $*$-algebra $\C[G_q/G_q^{\bar\lambda}]$ is generated by the elements $C^\lambda_{\eta,\zeta_\lambda}$ ($\eta\in V_\lambda$).
\end{enumerate}
\end{proposition}

Here by $G^\lambda_q$ we mean the quantization of $G^\lambda=K^{S,L}$ defined in Section~\ref{ssec:DJ}.

\bp
(1) It is well-known and not difficult to see that if we consider $V_\lambda$ as a module over the complexification $G_\C$ of $G$, then the subgroup $H$ of elements of $G_\C$ that leave the line $\C\xi_\lambda$ invariant is a parabolic subgroup, so $H\cap G=K^S$ for some $S\subset\Pi$. It is again easy to see using representation theory of $SU(2)$ that $S=\{\alpha\in\Pi: (\lambda,\alpha)=0\}$ and $\tilde K^S\subset G^\lambda$. Therefore $G^\lambda$ sits between $\tilde K^S$ and $K^S$, hence it coincides with $K^{S,L}$ for some subgroup $L\subset P(S^c)$. This means that $G^\lambda\cap T\subset T$ is the annihilator of $L\subset P$. Since by definition $G^\lambda\cap T$ is the kernel of $\lambda$ viewed as a character of $T$, its annihilator is $\Z\lambda$. Thus, $L=\Z\lambda$.

\smallskip

(2),(3) Since the involution maps the space spanned by the matrix coefficients $C^\lambda_{\eta,\xi_\lambda}$ ($\eta\in V_\lambda$) onto the space spanned by $C^{\bar\lambda}_{\eta,\zeta_{\bar\lambda}}$ ($\eta\in V_{\bar\lambda}$), part (2) for $\lambda$ is equivalent to part (3) for $\bar\lambda$. We also have an isomorphism $G_q\cong G_{q^{-1}}$ that swaps the meaning of highest and lowest weights. It follows that in order to prove (2) and (3) it suffices to establish (2) for $0<q\le 1$.

It is immediate that the elements $C^\lambda_{\eta,\xi_\lambda}$ lie in $\C[G_q/G_q^\lambda]$. Consider the $*$-algebra $\mathcal A_q$ generated by them. Note that it is unital by the unitarity of the representation of $G_q$ on $V_\lambda$. Denote by~$A_q$ the norm closure of $\mathcal A_q$ in $C(G_q)$. Since the algebra $\mathcal A_q$ is $G_q$-invariant with respect to the left action of $G_q$, the same is true for $A_q$ and, moreover, $\mathcal A_q$ coincides with the subalgebra of regular elements in $A_q$. It follows that it suffices to show that $A_q=C(G_q/G^\lambda_q)$.

For $q=1$ the equality $A_1=C(G/G^\lambda)$ follows from the Stone--Weierstrass theorem, since the functions $(\cdot\, \xi_\lambda,\eta)$ ($\eta\in V_\lambda$) separate points of $G/G^\lambda$.

For $0<q<1$ the argument is in principle similar, but now we have to use the Stone--Weierstrass theorem for type I C$^*$-algebras~\cite{MR0458185}*{Proposition~11.1.6}. For this we need to show that the irreducible representations of $C(G_q/G_q^\lambda)$ restrict to irreducible representations of $A_q$ and nonequivalent irreducible representations restrict to nonequivalent ones. The irreducible representations of $C(G_q/G_q^\lambda)=C(G_q/K^{S,L}_q)$ are described in~\cite{MR2914062}*{Theorem~2.2}, see also~\cite{MR1697598}*{Theorem 5.9}. Their equivalence classes are parameterized by the pairs $(w,x)\in W^S\times T/T_L$, where $W^S$ is the subset of the Weyl group consisting of elements $w$ such that $w\alpha>0$ for all $\alpha\in S$. Namely, as a representation corresponding to $(w,tT_L)$ one can take the restriction of Soibelman's representation $\pi_{w,t}$ of $C(G_q)$ to $C(G_q/G_q^\lambda)$. One can then check that the arguments used in the classification of irreducible representations of $C(G_q/G_q^\lambda)$ show that we have the required properties of $A_q$. Briefly, this goes as follows.

The proof of~\cite{MR1697598}*{Proposition 5.5}, which verifies irreducibility of the restriction of $\pi_{w,t}$ ($w\in W^S$, $t\in T$) to certain $*$-subalgebra of $C(G_q/G_q^\lambda)=C(G_q/K^{S,L}_q)$, uses only the matrix coefficients $C_{\eta,\xi_\lambda}^\lambda$ ($\eta\in V_\lambda$) and therefore shows that $\pi_{w,t}|_{A_q}$ is irreducible. In fact, it shows the stronger property that the restriction of $\pi_{w,t}$ to the C$^*$-subalgebra $\tilde A_q\subset A_q$ generated by the elements $(C_{\eta,\xi_\lambda}^\lambda)^*C_{\zeta,\xi_\lambda}^\lambda$ ($\eta,\zeta\in V_\lambda$) is irreducible.  Next, the proof of~\cite{MR1697598}*{Lemma 5.8} shows that the restrictions of $\pi_{w,t}$ and $\pi_{w',t'}$ to $A_q$ are inequivalent if $w,w'\in W^S$ and $w\ne w'$. Finally, take $w\in W^S$ and $t,t'\in T$. Then the restrictions of $\pi_{w,t}$ and $\pi_{w,t'}$ to $\tilde A_q$ coincide and, as we already observed, are irreducible. It follows that if the restrictions of $\pi_{w,t}$ and $\pi_{w,t'}$ to $A_q$ are equivalent, then they must be equal. Let $\eta\in V_\lambda$ be a unique up to a phase factor unit vector of weight $w\lambda$. Then it is known that $\pi_{w,t}(C^\lambda_{\eta,\xi_\lambda})\ne0$. As $\pi_{w,t'}(C_{\eta,\xi_\lambda}^\lambda)=\lambda(t^{-1}t')\pi_{w,t}(C_{\eta,\xi_\lambda}^\lambda)$, it follows that $t^{-1}t'\in \ker\lambda=T_L$. This completes the proof of the proposition.
\ep

Since the antipode transforms the matrix coefficients of $V_\lambda$ into those of $V_{\bar\lambda}$ and swaps the roles of highest and lowest weight vectors, we also get the following result.

\begin{corollary}\label{cor:stab}
For every $\lambda\in P_+$ and $q>0$, we have:
\begin{enumerate}
    \item[(1)] the $*$-algebra $\C[G_q^{\lambda}\backslash G_q]$ is generated by the elements $C^\lambda_{\xi_\lambda,\eta}$ ($\eta\in V_\lambda$);
    \item[(2)] the $*$-algebra $\C[G_q^{\bar\lambda}\backslash G_q]$ is generated by the elements $C^\lambda_{\zeta_\lambda,\eta}$ ($\eta\in V_\lambda$).
\end{enumerate}
\end{corollary}

\subsection{Commutation relations}\label{ssec:comm}
Fix $\lambda\in P_+$. We view $(V_{n\lambda})_{n\ge0}$ as a standard subproduct system, where the embedding $V_{n\lambda}\to V_\lambda^{\otimes n}$ is obtained by identifying $\xi_{n\lambda}$ with $\xi_\lambda^{\otimes n}$. Denote by $f_n$ the projection $V_\lambda^{\otimes n}\to V_{n\lambda}$. We will often decorate with the subscripts $\lambda$ and $q$ various associated constructions, so we write $\F_{\lambda,q}$, $\TT_{\lambda,q}$, $\OO_{\lambda,q}$, etc., for the corresponding Fock space, Toeplitz algebra, Cuntz--Pimsner algebra. For $q=1$ the subscript is omitted.

\begin{lemma}\label{lem:comm}
Take weight vectors $\xi,\zeta\in V_\lambda$. If either $\xi$ is a highest weight vector or $\zeta$ is a lowest weight vector, then in $\TT_{\lambda,q}$ we have
$$
S_\zeta S_\xi=q^{-(\lambda,\lambda)+(\wt(\zeta),\wt(\xi))}S_\xi S_\zeta.
$$
\end{lemma}

\bp
Consider the braiding $\sigma=\Sigma\RR$, where $\Sigma$ is the flip map. As $\sigma(\xi_\lambda\otimes\xi_\lambda)=q^{(\lambda,\lambda)}\xi_\lambda\otimes\xi_\lambda$ by~\eqref{eq:ermat2}, we have $f_2\sigma=q^{(\lambda,\lambda)}f_2$  on $V_\lambda\otimes V_\lambda$. Consider the map $S^{(2)}\colon V_\lambda\otimes V_\lambda\to\TT_{\lambda,q}$ defined by $S^{(2)}(\zeta\otimes\xi)=S_\zeta S_\xi$. We have $S^{(2)}f_2=S^{(2)}$ as a consequence of the inclusion $V_{n\lambda}\subset V_{2\lambda}\otimes V_{(n-2)\lambda}$. Hence $S^{(2)}\sigma=q^{(\lambda,\lambda)}S^{(2)}$. In other words, for all $\xi,\zeta\in V_\lambda$ we have
\begin{equation}\label{eq:comm}
S_\zeta S_\xi=q^{-(\lambda,\lambda)}S^{(2)}\sigma(\zeta\otimes\xi).
\end{equation}
Using again~\eqref{eq:ermat2} we get the result.
\ep

Using this lemma we can construct certain characters on $\TT_{\lambda,q}$, which will play an important role later.

\begin{lemma}\label{lem:char}
There exist states $\varphi_h$ and $\varphi_l$ on $\TT_{\lambda,q}$ such that
$$
\varphi_h(S_{\xi_\lambda})=1,\qquad \varphi_l(S_{\zeta_\lambda})=1.
$$
Moreover, any such states vanish on $\K(\F_{\lambda,q})$ and
\begin{enumerate}
    \item[(1)] if $q>1$, then $\varphi_h$ is a character satisfying $\varphi_h(S_\zeta)=0$ for all $\zeta\perp\xi_\lambda$;
    \item[(2)] if $q<1$, then $\varphi_l$ is a character satisfying $\varphi_l(S_\xi)=0$ for all $\xi\perp\zeta_\lambda$;
    \item[(3)] if $q=1$ and $\OO_\lambda$ is commutative, then $\varphi_h$ and $\varphi_l$ are characters satisfying $\varphi_h(S_\zeta)=0$ and $\varphi_l(S_\xi)=0$ for all $\zeta\perp\xi_\lambda$ and $\xi\perp\zeta_\lambda$.
\end{enumerate}
\end{lemma}

\bp
In order to show that $\varphi_h$ exists, consider the closed linear span $H\subset\F_{\lambda,q}$ of the vectors~$\xi_{n\lambda}$ ($n\ge0$). Consider the isometry $u\colon H\to H$ defined by $u\xi_{n\lambda}=\xi_{(n+1)\lambda}$. The C$^*$-algebra $C^*(u)\subset \B(H)$ generated by $u$ is the usual Toeplitz algebra $\TT$. Consider the state $\varphi$ on it such that $\varphi(u)=1$, namely, under the isomorphism $\TT/\K\cong C(\T)$ it corresponds to the evaluation at $1\in\T$. Extend it to a state on $\B(H)$. Denote by $P$ the projection $\F_{\lambda,q}\to H$. Then define a state~$\varphi_h$ on~$\TT_{\lambda,q}$ by $\varphi_h(a)=\varphi(PaP)$. As $PS_{\xi_\lambda}P=u$, it has the required property $\varphi_h(S_{\xi_\lambda})=1$. The existence of~$\varphi_l$ is proved in the same way by considering lowest weight vectors.

Assume now that $\varphi_h$ is any state such that  $\varphi_h(S_{\xi_\lambda})=1$. Since $S_{\xi_\lambda}$ is a contraction, it follows that $\varphi_h(S_{\xi_\lambda}^*S_{\xi_\lambda})=1=\varphi_h(S_{\xi_\lambda}S_{\xi_\lambda}^*)$ and $S_{\xi_\lambda}$ lies in the multiplicative domain of $\varphi_h$.

By~\eqref{eq:e-0} we have $\sum_i S_iS_i^*=1-e_0$ for any choice of an orthonormal basis in $V_\lambda$. Choosing a basis containing $\xi_\lambda$ and applying $\varphi_h$ we conclude that $\varphi_h(e_0)=0$ and $\varphi_h(S_\zeta S^*_\zeta)=0$ for any $\zeta\perp\xi_\lambda$. More generally, denote by $e_n$ the projection $\F_{\lambda,q}\to V_{n\lambda}$. From~\eqref{eq:e-0} we get by induction that
$$
\sum_{i_1,\dots,i_n}S_{i_1}\dots S_{i_n}S_{i_n}^*\dots S_{i_1}^*=1-e_0-\dots-e_{n-1},
$$
and since $\varphi_h(S_{\xi_\lambda}^n(S_{\xi_\lambda}^*)^n)=1$, we conclude that $\varphi_h(e_k)=0$ for all $0\le k<n$. Hence $\varphi_h$ vanishes on~$\K(\F_{\lambda,q})$. In the same way one checks that $\varphi_l$ vanishes on~$\K(\F_{\lambda,q})$.

Next, assume $q>1$. Since we already know that  $S_{\xi_\lambda}$ lies in the multiplicative domain of $\varphi_h$ and $\varphi_h(S_\zeta S^*_\zeta)=0$ (and hence $\varphi_h(S_\zeta)=0$) for any $\zeta\perp\xi_\lambda$, in order to prove (1) it remains to show that every such $S_\zeta$ lies in the multiplicative domain of $\varphi_h$, that is, $\varphi_h(S_\zeta^* S_\zeta)=0$. We may assume that $\zeta\in V_\lambda(\mu)\ne0$ for some $\mu\ne\lambda$. Then using Lemma~\ref{lem:comm} we get
$$
\varphi_h(S_\zeta^* S_\zeta)=\varphi_h(S_{\xi_\lambda}^*S_\zeta^* S_\zeta S_{\xi_\lambda})=q^{-2(\lambda,\lambda)+2(\mu,\lambda)}\varphi_h(S_\zeta^*S_{\xi_\lambda}^*S_{\xi_\lambda} S_\zeta )
\le q^{-2(\lambda,\lambda)+2(\mu,\lambda)}\varphi_h(S_\zeta^*S_\zeta ).
$$
As $q>1$ and $(\lambda-\mu,\lambda)>0$, this implies that $\varphi_h(S_\zeta^* S_\zeta)=0$.

Part (2) is proved similarly, using that $\wt(\zeta_\lambda)=w_0\lambda$ and
$$
(\lambda,\lambda)-(\mu,w_0\lambda)=(\lambda-w_0\mu,\lambda)>0
$$
for all weights $\mu\ne w_0\lambda$ in the weight decomposition of $V_\lambda$.

Finally, in order to prove (3) for $\varphi_h$ it suffices to show $\varphi_h(S_\zeta^* S_\zeta)=0$ for all $\zeta\perp\xi_\lambda$, as we observed above. But this is clearly true, since we already know that $\varphi_h(S_\zeta S^*_\zeta)=0$, $\varphi_h$ factors through~$\OO_\lambda$ and $\OO_\lambda$ is commutative by assumption. Similarly for $\varphi_l$.
\ep

\subsection{Approximate commutation relations for creation and annihilation operators}\label{ssec:star-comm}
In this subsection it will be convenient to consider a setup that is more general than we strictly speaking need for our analysis of Cartan subproduct systems.

For a fixed $q>0$, we view the collection $(V_\lambda)_{\lambda\in P_+}$ as a standard subproduct system over the monoid $P_+$, where the embedding $V_{\lambda+\mu}\to V_\lambda\otimes V_\mu$ is defined by identifying $\xi_{\lambda+\mu}$ with $\xi_\lambda\otimes\xi_\mu$. Similarly to subproduct systems over $\Z_+$, we then get a full Fock space
$$
\F_{P_+,q}=\bigoplus_{\lambda\in P_+}V_\lambda
$$
and creation operators $L_\xi\colon \F_{P_+,q}\to \F_{P_+,q}$ mapping $V_\mu$ into $V_{\lambda+\mu}$ for $\xi\in V_\lambda$.

Denote by $P^h_{\lambda,\mu}\colon V_\lambda\otimes V_\mu\to V_\lambda\otimes V_\mu$ the orthogonal projection onto the subspace of vectors killed by all the $E_i$'s. Denote by $P^h_\mu\colon V_\mu\to V_\mu$ the orthogonal projection onto $\mathbb C\xi_\mu$; in other words, $P^h_\mu=P^h_{0,\mu}$. Similarly, denote by $P^l_{\lambda,\mu}\colon V_\lambda\otimes V_\mu\to V_\lambda\otimes V_\mu$ the orthogonal projection onto the subspace of vectors killed by all the $F_i$'s, and let $P^l_\mu=P^l_{0,\mu}$. Denote by $f_{\lambda,\mu}\colon V_\lambda\otimes V_\mu\to V_\lambda\otimes V_\mu$ the projection onto the Cartan component of the tensor product.

\begin{proposition} \label{prop:star-commute}
Take $q>0$ and $\lambda\in P_+$. For $\mu\in\lambda+P_+$, define linear maps $A_\mu\colon V_\lambda\otimes \bar V_\lambda\to \B(V_\mu)$ and $B_\mu\colon \bar V_\lambda\otimes V_\lambda\to \B(V_\mu)$ by
$$
A_\mu(\xi\otimes\bar\zeta)=L_\xi L_\zeta^*|_{V_\mu},\qquad B_\mu(\bar\zeta\otimes\xi)=L_\zeta^* L_\xi|_{V_\mu}.
$$
Then there exists a constant $C>0$ depending only on $q$ and $\lambda$ such that, for all $\mu\in\lambda+P_+$,
\begin{align*}
\|B_\mu-q^{-(\lambda,\lambda)}A_\mu\sigma^{-1}\|&\le C(\big\|f_{\lambda,\mu}\big((1-P^h_\lambda)\otimes P^h_\mu\big)\big\|+\big\|f_{\lambda,\mu-\lambda}\big((1-P^h_\lambda)\otimes P^h_{\mu-\lambda}\big)\big\|\\
& \qquad+\|(1-1\otimes P^h_\mu)P^h_{\lambda,\mu}\|\big),\\
\|B_\mu-q^{(\lambda,\lambda)}A_\mu\sigma\|&\le C(\big\|f_{\lambda,\mu}\big((1-P^l_\lambda)\otimes P^l_\mu\big)\big\|+\big\|f_{\lambda,\mu-\lambda}\big((1-P^l_\lambda)\otimes P^l_{\mu-\lambda}\big)\big\|\\
& \qquad+\|(1-1\otimes P^l_\mu)P^l_{\lambda,\mu}\|\big),
\end{align*}
where $\sigma=\Sigma\RR$ is the braiding on $\Rep G_q$ and the spaces $\bar V_\lambda\otimes V_\lambda$ and $\B(V_\mu)$ are considered with their usual Hilbert space and operator norms, resp.
\end{proposition}

Here $\bar V_\lambda$ is equipped with the contragredient $U_q\mathfrak g$-action, that is, $\omega\bar\xi=\overline{\hat S_q(\omega)^*\xi}$.

\bp
We will only prove the first inequality, the second one is proved similarly by considering the lowest weight vectors instead of the highest weight vectors.

Let us show first that for any weight vectors $\xi,\zeta\in V_\lambda$ of norm one there exists a constant $C_1>0$ such that for all $\mu\in\lambda+P_+$ we have
\begin{multline}\label{eq:convergence1}
\|(B_\mu\sigma)(\xi\otimes\bar\zeta)\xi_\mu-q^{-(\lambda,\lambda)}A_\mu(\xi\otimes\bar\zeta)\xi_\mu\|\le C_1\big(\big\|f_{\lambda,\mu}\big((1-P^h_\lambda)\otimes P^h_\mu\big)\big\|\\ +\big\|f_{\lambda,\mu-\lambda}\big((1-P^h_\lambda)\otimes P^h_{\mu-\lambda}\big)\big\|\big).
\end{multline}

From~\eqref{eq:ermat1} we know that $\sigma(\xi\otimes\bar\zeta)$ is the sum of $q^{(\wt(\xi),\wt(\bar\zeta))}\bar\zeta\otimes\xi=q^{-(\wt(\xi),\wt(\zeta))}\bar\zeta\otimes\xi$
and vectors of the form $\bar{\zeta'}\otimes\xi'$ such that $\zeta'$ and $\xi'$ have weights strictly lower than~$\wt(\zeta)$ and~$\wt(\xi)$, resp. It follows that in order to prove~\eqref{eq:convergence1} it suffices to show that if $\wt(\xi)=\wt(\zeta)=\lambda$, then
\begin{equation}\label{eq:convergence2-0}
B_\mu(\bar\zeta\otimes\xi)\xi_\mu=A_\mu(\xi\otimes\bar\zeta)\xi_\mu,
\end{equation}
while if either $\wt(\xi)$ or $\wt(\zeta)$ is different from $\lambda$, then
\begin{align}
\|B_\mu(\bar\zeta\otimes\xi)\xi_\mu\|&\le \big\|f_{\lambda,\mu}\big((1-P^h_\lambda)\otimes P^h_\mu\big)\big\|,\label{eq:convergence2-b}\\
\|A_\mu(\xi\otimes\bar\zeta)\xi_\mu\|&\le \big\|f_{\lambda,\mu-\lambda}\big((1-P^h_\lambda)\otimes P^h_{\mu-\lambda}\big)\big\|\label{eq:convergence2-a}
\end{align}
for all $\mu\in\lambda+P_+$.

In order to check~\eqref{eq:convergence2-0} we may assume that $\xi=\zeta=\xi_\lambda$. Then both sides of~\eqref{eq:convergence2-0} equal~$\xi_\mu$.

Next, let us prove~\eqref{eq:convergence2-b}, inequality~\eqref{eq:convergence2-a} is proved similarly. Assume first that $\nu=\wt(\xi)$ is different from~$\lambda$. Then
$$
\|B_\mu(\bar\zeta\otimes\xi)\xi_\mu\|\le\|L_\xi\xi_\mu\|=\|f_{\lambda,\mu}(\xi\otimes\xi_\mu)\|\le\big\|f_{\lambda,\mu}\big((1-P^h_\lambda)\otimes P^h_\mu\big)\big\|.
$$
Assume now that $\wt(\xi)=\lambda$, but $\wt(\zeta)\ne\lambda$. We may assume then that $\xi=\xi_\lambda$, so that $L_\xi\xi_\mu=\xi_{\lambda+\mu}$. Then $B_\mu(\bar\zeta\otimes\xi)\xi_\mu=L_\zeta^* \xi_{\lambda+\mu}=0$, since $\xi_{\lambda+\mu}=\xi_\lambda\otimes\xi_\mu$ is orthogonal to $\mathbb C\zeta\otimes V_\mu$. This finishes the proof of~\eqref{eq:convergence2-b} and~\eqref{eq:convergence2-a}, hence also the proof of~\eqref{eq:convergence1}.

\smallskip

Observe next that the maps $A_\mu$ and $B_\mu$ become $G_q$-equivariant if $\B(V_\mu)$ is equipped with the adjoint action, so that $(\operatorname{ad}\omega)(T)=\pi_\mu(\omega_{(1)})T \pi_\mu(\hat S_q(\omega_{(2)}))$ for $\omega\in U_q\g$ and $T\in \B(V_\mu)$, where~$\pi_\mu$ denotes the representation of $U_q\g$ on $V_\mu$. Choose an orthonormal basis $(\xi_i)_i$ in $V_\lambda$. Define linear maps $\tilde A_\mu, \tilde B_\mu\colon V_\lambda\otimes V_\mu\to V_\lambda\otimes V_\mu$ by
$$
\tilde A_\mu(\xi\otimes\zeta)=q^{-(\lambda,\lambda)}\sum_i\xi_i\otimes (A_\mu\sigma^{-1})(\bar\xi_i\otimes\xi)\zeta,\quad
\tilde B_\mu(\xi\otimes\zeta)=\sum_i\xi_i\otimes B_\mu(\bar\xi_i\otimes\xi)\zeta.
$$
Using that the vector $\sum_i\xi_i\otimes\bar\xi_i$ is invariant it is easy to check that $\tilde A_\mu$ and $\tilde B_\mu$ are morphisms of $U_q\g$-modules. Then~\eqref{eq:convergence1} can be rephrased by saying that there exists a constant $C_2>0$ depending only on $q$ and $\lambda$ such that
\begin{equation}\label{eq:convergence3}
\|(\tilde B_\mu-\tilde A_\mu)|_{V_\lambda\otimes\C\xi_\mu}\|\le C_2\big(\big\|f_{\lambda,\mu}\big((1-P^h_\lambda)\otimes P^h_\mu\big)\big\| +\big\|f_{\lambda,\mu-\lambda}\big((1-P^h_\lambda)\otimes P^h_{\mu-\lambda}\big)\big\|\big).
\end{equation}
What we need is to estimate the norm of $\tilde B_\mu-\tilde A_\mu$ on the entire space $V_\lambda\otimes V_\mu$.

Note that by definition there is a universal (depending only on $q$ and $\lambda$) bound $C_3$ on the norm of $\tilde B_\mu-\tilde A_\mu$. By decomposing the module $V_\lambda\otimes V_\mu$ into simple ones we can identify it with $\bigoplus_\nu V_\nu\otimes H_\nu$ for some Hilbert spaces $H_\nu$ with trivial $G_q$-action. Under this identification any morphism $T\colon V_\lambda\otimes V_\mu\to V_\lambda\otimes V_\mu$ has the form $\sum_\nu 1\otimes T_\nu$ for some $T_\nu\in \B(H_\nu)$, while $P^h_{\lambda,\mu}$ becomes the projection onto $\bigoplus_\nu \C\xi_\nu\otimes H_\nu$. It follows that $\|T\|=\|TP^h_{\lambda,\mu}\|$. Therefore
$$
\|\tilde B_\mu-\tilde A_\mu\|=\|(\tilde B_\mu-\tilde A_\mu)P^h_{\lambda,\mu}\|\le C_3\|(1-1\otimes P^h_\mu)P^h_{\lambda,\mu}\|+\|(\tilde B_\mu-\tilde A_\mu)(1\otimes P^h_\mu)\|.
$$
Combined with~\eqref{eq:convergence3} this finishes the proof of the proposition.
\ep

\subsection{Asymptotics of Clebsch--Gordan coefficients}\label{ssec:conj}
We conjecture that when $q\ge1$ and $n$ is large, then every highest weight vector in $V_\lambda\otimes V_{n\lambda}$ is close to a vector of the form $\xi\otimes \xi_{n\lambda}$, and any vector of this form with $\xi\perp\xi_\lambda$ is almost orthogonal to $V_{(n+1)\lambda}$. This can be formulated in terms of the Clebsch--Gordan coefficients of~$G_q$. For example, the first part means that, given an embedding $V_\nu\hookrightarrow V_\lambda\otimes V_{n\lambda}$ and an orthonormal basis $(\xi_i)_i$ in $V_\lambda$, we have $\sum_i|(\xi_\nu,\xi_i\otimes\xi_{n\lambda})|^2\approx1$. It is, however, more convenient to use the projections $P^h_{\lambda,\mu}$, $P^l_{\lambda,\mu}$ and $f_{\lambda,\mu}$ introduced before Proposition~\ref{prop:star-commute}. Here is the precise formulation of our conjecture for $q>0$ and $\lambda\in P_+$.

\begin{conjecture}\label{conj}
The following properties hold: if $q\ge1$, then
\begin{equation}\label{eq:conj1}
\|(1-1\otimes P^h_{n\lambda})P^h_{\lambda,n\lambda}\|\to0\quad\text{and}\quad \big\|f_{\lambda,n\lambda}\big((1-P^h_\lambda)\otimes P^h_{n\lambda}\big)\big\|\to0\quad\text{as}\quad n\to+\infty,
\end{equation}
and if $q\le1$, then
\begin{equation}\label{eq:conj2}
\|(1-1\otimes P^l_{n\lambda})P^l_{\lambda,n\lambda}\|\to0\quad\text{and}\quad \big\|f_{\lambda,n\lambda}\big((1-P^l_\lambda)\otimes P^l_{n\lambda}\big)\big\|\to0\quad\text{as}\quad n\to+\infty. 
\end{equation}
Moreover, for $q\ne1$, the convergences are exponentially fast, that is, for $q>1$ we have
\begin{equation}\label{eq:conjS}
\|(1-1\otimes P^h_{n\lambda})P^h_{\lambda,n\lambda}\|\le Ct^n\quad\text{and}\quad\big\|f_{\lambda,n\lambda}\big((1-P^h_\lambda)\otimes P^h_{n\lambda}\big)\big\|\le Ct^n
\end{equation}
for some $C>0$ and $0<t<1$, and similarly for $q<1$.
\end{conjecture}

Properties \eqref{eq:conj1} for $q\ge 1$ and \eqref{eq:conj2} for $q^{-1}$ are, in fact, equivalent thanks to the standard isomorphism $G_q\cong G_{q^{-1}}$. Using the unitary antipode it also possible to formulate the conjecture in terms of the modules $V_{n\lambda}\otimes V_\lambda$. It should be stressed that a similar conjecture with the roles of highest and lowest weights swapped is \emph{not} true for $q\ne1$ already for $G=SU(2)$. Therefore a proper formulation of the conjecture depends on the conventions for the $q$-deformations.

Recall that a weight $\lambda\in P_+$ is called regular if $\lambda(i)>0$ for all $i$. Denote by $P_{++}\subset P_+$ the subset of regular weights. For regular weights the conjecture can be formulated in a neater way. We will need the following elementary lemma.

\begin{lemma}
Assume $H$ is a Hilbert space and $p,q\in B(H)$ are projections of the same finite rank. Then
$$
\|(1-p)q\|\le\|p-q\|\le 2\|(1-p)q\|.
$$
\end{lemma}

\bp
The first inequality holds for all projections, since $(1-p)q=(q-p)q$. In order to prove the second inequality, observe that if we view $pq$ as an operator $T\colon qH\to pH$ between finite dimensional spaces of the same dimension, then $T^*T$ and $TT^*$ have equal spectra, hence $\|1-T^*T\|=\|1-TT^*\|$, that is, $\|q(1-p)q\|=\|p(1-q)p\|$. It follows that
$$
\|p(1-q)\|=\|p(1-q)p\|^{1/2}=|q(1-p)q\|^{1/2}=\|(1-p)q\|,
$$
and then $\|p-q\|\le\|p(1-q)\|+\|(1-p)q\|=2\|(1-p)q\|$.
\ep

Now, if $\lambda\in P_{++}$, then the first convergence in~\eqref{eq:conj1} implies the second one, and the conjecture can be formulated as follows.

\begin{lemma}
Assume  $q\ge1$ and $\lambda\in P_{++}$. Then~\eqref{eq:conj1} is equivalent to the following property: $\|P^h_{\lambda,n\lambda}-1\otimes P^h_{n\lambda}\|\to0$ as $n\to+\infty$. Moreover, for $q>1$ the last convergence is exponentially fast if only if \eqref{eq:conjS} is true.
\end{lemma}

\bp
By Theorem~\ref{thm:mult}, when $n$ is sufficiently large, the projections $P^h_{\lambda,n\lambda}$ and $1\otimes P^h_{n\lambda}$ have the same rank $\dim V_\lambda$. From this and the previous lemma we then see that the convergences $\|(1-1\otimes P^h_{n\lambda})P^h_{\lambda,n\lambda}\|\to0$ and $\|P^h_{\lambda,n\lambda}-1\otimes P^h_{n\lambda}\|\to0$ are equivalent, and if one is exponentially fast, then the other is exponentially fast as well. Next, we have
$$
\big\|f_{\lambda,n\lambda}\big((1-P^h_\lambda)\otimes P^h_{n\lambda}\big)\big\|\le \|f_{\lambda,n\lambda}(P^h_{\lambda,n\lambda}-P^h_\lambda\otimes P^h_{n\lambda})\|+\|P^h_{\lambda,{n\lambda}}-1\otimes P^h_{n\lambda}\|.
$$
The operator $P^h_{\lambda,n\lambda}-P^h_\lambda\otimes P^h_{n\lambda}$ is the projection onto the space spanned by the highest weight vectors in $V_\lambda\otimes V_{n\lambda}$ orthogonal to $\xi_\lambda\otimes\xi_{n\lambda}$, that is, the highest weight vectors of all components of the tensor product that are different from the Cartan component. Therefore $f_{\lambda,n\lambda}(P^h_{\lambda,n\lambda}-P^h_\lambda\otimes P^h_{n\lambda})=0$, so that $\big\|f_{\lambda,n\lambda}\big((1-P^h_\lambda)\otimes P^h_{n\lambda}\big)\big\|\le \|P^h_{\lambda,n\lambda}-1\otimes P^h_\mu\|$ and the second convergence in~\eqref{eq:conj1} is implied by the first one.
\ep

In order to prove the conjecture for regular weights it is enough to verify a stronger property for just one (in general, not irreducible) representation. Namely, given a finite dimensional unitary $G_q$-module $V$, denote by $P^h_{V,\mu}$ the projection onto the subspace of vectors in $V\otimes V_\mu$ killed by all the $E_i$'s. Let us say that a $G_q$-module $V$ is faithful if every module~$V_\lambda$ is contained in some tensor power of $V$. It is not difficult to see that $V$ is faithful if and only if the highest weights appearing in the decomposition of $V$ into simple modules generate $P$ modulo the root lattice $Q$. For $\mu\in P_+
$, denote by $N(\mu)$ the minimum of $\mu(i)$ over all $1\le i\le r$.

\begin{lemma}\label{lem:regular}
Assume $q\ge1$ and $V$ is a faithful finite dimensional unitary $G_q$-module such that $\|P^h_{V,\mu}-1\otimes P^h_\mu\|\to0$ as $N(\mu)\to+\infty$, and the convergence is exponentially fast when $q>1$, meaning that $\|P^h_{V,\mu}-1\otimes P^h_\mu\|\le Ct^{N(\mu)}$ for some $C>0$ and $0<t<1$. Then Conjecture~\ref{conj} is true for $q^{\pm1}$ and all $\lambda\in P_{++}$.
\end{lemma}

\bp
For a finite dimensional unitary $G_q$-module $W$ consider the following property:
\begin{equation}\label{eq:conjW}
\|P^h_{W,\mu}-1\otimes P^h_\mu\|\to0\quad\text{as}\quad N(\mu)\to+\infty.
\end{equation}
By assumption this holds for $W=V$. Moreover, we have
$$
\|P^h_{V\otimes W,\mu}-1\otimes P^h_{W,\mu}\|\to0\quad\text{as}\quad N(\mu)\to+\infty
$$
for any $W$, since by Theorem~\ref{thm:mult} the module $W\otimes V_\mu$ decomposes into a finite sum, with a universal (depending only on $W$) bound  on the number of summands, of simple modules with highest weights $\nu$ such that the differences $N(\nu)-N(\mu)$ are bounded. It follows that if~\eqref{eq:conjW} holds for some $W$, then it holds for $V\otimes W$ as well. By induction we then conclude that~\eqref{eq:conjW} holds for $W=V^{\otimes n}$ for all $n$. Since by assumption any $V_\lambda$ is contained in $V^{\otimes n}$ for some $n$, it follows that~\eqref{eq:conjW} holds for $W=V_\lambda$ for all $\lambda\in P_+$, that is, $\|P^h_{\lambda,\mu}-1\otimes P^h_\mu\|\to0$ as $N(\mu)\to+\infty$. The same argument shows that if $\|P^h_{V,\mu}-1\otimes P^h_\mu\|\to0$ exponentially fast, then the convergence $\|P^h_{\lambda,\mu}-1\otimes P^h_\mu\|\to0$ as $N(\mu)\to+\infty$ is exponentially fast as well for all $\lambda\in P_+$. In particular, Conjecture~\ref{conj} is true for $q$ and all $\lambda\in P_{++}$.
\ep

Using this lemma and known formulas for Clebsch--Gordan coefficients of certain tensor product modules, one can verify the conjecture in some cases. In particular, we have the following result.

\begin{theorem}\label{thm:SUN}
Conjecture~\ref{conj} is true for $G=SU(N)$, $q>0$, $\lambda\in P_{++}$ and $\lambda\in\mathbb N\omega_1$.
\end{theorem}

Here $\omega_1$ denotes the first fundamental weight of $SU(N)$ with respect to the standard choice of simple roots. We defer the proof of the theorem to Appendix~\ref{appendix}.

\subsection{Cuntz--Pimsner algebras of Cartan subproduct systems}
We are now ready to prove our main result.

\begin{theorem}\label{thm:main}
Assume $G$ is a simply connected semisimple compact Lie group, $q>0$ and $\lambda\in P_+\setminus\{0\}$. Consider the corresponding Cartan subproduct system $(V_{n\lambda})_{n\ge0}$ of $G_q$-modules and the associated Cuntz--Pimsner algebra $\OO_{\lambda,q}$. Then
\begin{enumerate}
    \item[(1)] if $q\ge1$ and~\eqref{eq:conj1} holds, then we have a $G_q$-equivariant isomorphism $\OO_{\lambda,q}\cong C(G^{\lambda}_q\backslash G_q)$, $s_\eta\mapsto C^\lambda_{\xi_\lambda,\eta}$ ($\eta\in V_\lambda$);
    \item[(2)] if $q\le1$ and~\eqref{eq:conj2} holds, then we have a $G_q$-equivariant isomorphism $\OO_{\lambda,q}\cong C(G^{\bar\lambda}_q\backslash G_q)$, $s_\eta\mapsto C^\lambda_{\zeta_\lambda,\eta}$ ($\eta\in V_\lambda$).
\end{enumerate}
\end{theorem}

\bp
We will prove (1), part (2) is proved similarly. Consider the linear maps $A\colon V_\lambda\otimes\bar V_\lambda\to\TT_{\lambda,q}$ and $B\colon \bar V_\lambda\otimes V_\lambda\to\TT_{\lambda,q}$ defined by
$A(\xi\otimes\bar\zeta)=S_\xi S^*_\zeta$ and $B(\bar\zeta\otimes\xi)=S_\zeta^*S_\xi$. Condition~\eqref{eq:conj1} and Proposition~\ref{prop:star-commute} imply that for all $\xi,\zeta\in V_\lambda$ we have
$$
B(\bar\zeta\otimes\xi)=q^{-(\lambda,\lambda)}A\sigma^{-1}(\bar\zeta\otimes\xi)\quad\operatorname{mod}\quad \K(\F_{\lambda,q}).
$$
It follows that $\OO_{\lambda,q}$ admits normal ordering. Moreover, when $q=1$ and therefore $\sigma$ is just the flip map, this shows that the elements~$s_\zeta^*$ and~$s_\xi$ commute, and since the elements $s_\zeta$ and $s_\xi$ commute as well by~\eqref{eq:comm}, we conclude that $\OO_\lambda$ is a commutative C$^*$-algebra.

By Lemma~\ref{lem:char}, there is a character $\varphi_h$ on $\OO_{\lambda,q}$ such that $\varphi_h(s_{\xi_\lambda})=1$ and $\varphi_h(s_\eta)=0$ for all $\eta\perp\xi_\lambda$. Denoting by $\alpha\colon\OO_{\lambda,q}\to\OO_{\lambda,q}\otimes C(G_q)$ the action of $G_q$, define a unital $*$-homomorphism
$$
\pi\colon\OO_{\lambda,q}\to C(G_q)\quad\text{by}\quad \pi(a)=(\varphi_h\otimes\iota)\alpha(a).
$$
Then we get $\pi(s_\eta)=C^\lambda_{\xi_\lambda,\eta}$ for all $\eta\in V_\lambda$. By Corollary~\ref{cor:stab}(1) we conclude that the image of $\pi$ equals $C(G^\lambda_q\backslash G_q)$.

The homomorphism $\pi\colon \OO_{\lambda,q}\to C(G^\lambda_q\backslash G_q)$ is $G_q$-equivariant. By Proposition~\ref{prop:normal}, the action of $G_q$ on $\OO_{\lambda,q}$ is ergodic. Since $G_q$ is coamenable by~\cite{neshveyev-tuset-book}*{Theorem~2.7.14}, a standard argument shows then that $\pi$ is injective, see, e.g., the discussion preceding Definition~1.2 in \cite{MR4919591}*{Section~1}. This completes the proof of the theorem.
\ep

By Theorem~\ref{thm:SUN} this result applies unconditionally to $G=SU(N)$, $q>0$, $\lambda\in P_{++}$ and $\lambda\in\mathbb N\omega_1$. For $\lambda=\omega_1$ this gives again Theorem~\ref{thm:arveson}, but the two proofs are not independent, since some of the computations in Section~\ref{sec:q-sym} are used in the proof of Theorem~\ref{thm:SUN}.

\bigskip

\section{Gauge-invariant subalgebra}\label{sec:gauge}

Fix again a simply connected semisimple compact Lie group $G$ and a dominant integral weight $\lambda\in P_+\setminus\{0\}$. In this section we will study the gauge-invariant subalgebra of $\TT_{\lambda,q}$. We will assume that $q\ge1$, but as should be clear from the previous two sections, similar results can be obtained for $0<q\le1$ as well by swapping the roles of highest and lowest weights, or simply by using the isomorphism $G_{q^{-1}}\cong G_q$.

\subsection{A continuous field structure and Berezin quantization}

Consider the unitary representation $z\mapsto u_z$ of $\T$ on $\F_{\lambda,q}$, where $u_z$ is the unitary that acts on~$V_{n\lambda}$ as the multiplication by $z^n$. Then we get an action $\Ad u$ of $\T$ on $\B(\F_{\lambda,q})$. The automorphisms $\Ad u_z$ leave~$\TT_{\lambda,q}$ globally invariant and define the gauge action~$\gamma$ of $\T$ on $\TT_{\lambda,q}$, namely, $\gamma_z(S_\xi)=zS_\xi$ for $\xi\in V_\lambda$. This action passes to $\OO_{\lambda,q}$. We denote by $\TT^{(0)}_{\lambda,q}\subset\TT_{\lambda,q}$ and $\OO^{(0)}_{\lambda,q}\subset\OO_{\lambda,q}$ the gauge-invariant subalgebras.

We assume that condition \eqref{eq:conj1} is satisfied. Then by Theorem~\ref{thm:main} we have a surjective homomorphism $\pi\colon\TT_{\lambda,q}\to C(G^\lambda_q\backslash G_q)$ with kernel $\K(\F_{\lambda,q})$. Let $S=\{\alpha\in\Pi: (\lambda,\alpha)=0\}$.

\begin{lemma}\label{lem:flag}
The homomorphism $\pi\colon\TT_{\lambda,q}\to C(G^\lambda_q\backslash G_q)$ maps $\TT^{(0)}_{\lambda,q}\subset\TT_{\lambda,q}$ onto $C(K^S_q\backslash G_q)$ and therefore defines an isomorphism $\OO^{(0)}_{\lambda,q}\cong C(K^S_q\backslash G_q)$.
\end{lemma}

\bp
Our fixed maximal torus $T\subset G$ acts by left translations on $C(G_q)$. Explicitly, an element $t=\exp(X)\in T$ rescales a matrix coefficient $C^\mu_{\zeta,\xi}$, where $\zeta\in V_\mu$ is a weight vector, by $e^{-\wt(\zeta)(X)}$. In particular, each $t\in T$ rescales all elements $C^\lambda_{\xi_\lambda,\eta}=\pi(S_\eta)$, $\eta\in V_\lambda$, in the same way. It follows
that $\pi$ defines an isomorphism of $\OO^{(0)}_{\lambda,q}$ onto $C(G^\lambda_q\backslash G_q)\cap C(T\backslash G_q)$. Recalling that $G^\lambda_q=K^{S,L}_q$, where $L=\Z\lambda$, we see that the last intersection coincides with $C(K^S_q\backslash G_q)$, which proves the lemma.
\ep

The homomorphism $\TT^{(0)}_{\lambda,q}\to C(K^S_q\backslash G_q)$ is closely related to Berezin quantization (or rather dequantization). Namely, for every $n$ we can define a $G_q$-equivariant ucp map
$$
\sigma_n\colon \B(V_{n\lambda})\to C(K^S_q\backslash G_q)\quad\text{by}\quad \sigma_n(T)=(\omega_n\otimes\iota)(U_n(T\otimes1)U_n^*),
$$
where $\omega_n=(\cdot\,\xi_{n\lambda},\xi_{n\lambda})$ and $U_n\in \B(V_{n\lambda})\otimes \C[G_q]$ is the unitary defining the representation of~$G_q$ on~$V_{n\lambda}$. In the classical case $q=1$ the function $\sigma_n(T)$ on $G$ is called a covariant Berezin symbol of~$T$.

The subspaces $V_{n\lambda}\subset\F_{\lambda,q}$ are invariant under $\TT^{(0)}_{\lambda,q}$. Given $x\in\TT^{(0)}_{\lambda,q}$, let us write $x_n$ for $x|_{V_{n\lambda}}$.

\begin{lemma} \label{lem:norm-convergence}
For every $x\in\TT^{(0)}_{\lambda,q}$, the sequence $(\sigma_n(x_n))_n$ converges in norm to $\pi(x)$.
\end{lemma}

\bp
Recall that the homomorphism $\pi$ from Theorem~\ref{thm:main}(1) is defined by $\pi(x)=(\varphi_h\otimes\iota)\alpha(x)$. As we showed in the proof of Proposition~\ref{prop:star-commute}, if $\zeta\in V_\lambda$ is a weight vector of weight different from~$\lambda$, then $S_\zeta^*\xi_{n\lambda}=0$ and $\|S_\zeta\xi_{n\lambda}\|\le\|f_{\lambda,n\lambda}\big((1-P^h_\lambda)\otimes P^h_{n\lambda}\big)\|\,\|\zeta\|$. From this it is not difficult to see that the restriction of the character $\varphi_h$ on $\TT_{\lambda,q}$ to $\TT_{\lambda,q}^{(0)}$ is given by $\varphi_h(x)=\lim_n\omega_n(x_n)$. It follows that for every $x\in\TT_{\lambda,q}^{(0)}$ the sequence $(\sigma_n(x_n))_n$ converges weakly to $\pi(x)\in C(K^S_q\backslash G_q)$.
The convergence is in fact in norm, since every element of $\TT_{\lambda,q}^{(0)}$ can be approximated by elements $x$ contained in finite dimensional $G_q$-invariant subspaces and for such~$x$ the claim follows immediately by the definition of $\sigma_n$ and $\pi$.
\ep

Therefore $\pi|_{\TT^{(0)}_{\lambda,q}}$ can be seen as a limit of quantum analogues of covariant Berezin symbols. This has the following consequence.

\begin{theorem}\label{thm:cont-field}
Let $G$ be a simply connected semisimple compact Lie group, $\lambda\in P_+\setminus\{0\}$ and $q\ge1$. Let $S=\{\alpha\in\Pi: (\lambda,\alpha)=0\}$. Assume condition~\eqref{eq:conj1} is satisfied. Then $\TT_{\lambda,q}^{(0)}$
is the section algebra of a continuous field of C$^*$-algebras over $\bar\Z_+=\Z_+\cup\{\infty\}$ with fibers $\B(V_{n\lambda})$ at $n\in\Z_+$ and $C(K^S_q\backslash G_q)$ at $\infty$.
\end{theorem}

\bp
We view every $x\in \TT^{(0)}_{\lambda,q}$ as a section of the field of C$^*$-algebras over $\bar\Z_+$ taking value~$x_n$ at $n\in\Z_+$ and $x_\infty=\pi(x)$ at $\infty$. In order to show that this way we get a continuous field structure as in the statement of the theorem we need to check three properties: 1) the map $\TT^{(0)}_{\lambda,q}\ni x\mapsto x_t$ has image $\B(V_{n\lambda})$ for $t=n$ and $C(K^S_q\backslash G_q)$ for $t=\infty$; 2) the collection of sections defined by the elements of $\TT^{(0)}_{\lambda,q}$ is stable under multiplication by continuous functions on~$\bar\Z_+$; and 3) the map $t\mapsto\|x_t\|$ is continuous on $\bar\Z_+$ for all $x\in \TT^{(0)}_{\lambda,q}$.

The first two properties follow from Lemma~\ref{lem:flag} and the fact $\TT^{(0)}_{\lambda,q}$ contains all gauge-invariant compact operators on $\F_{\lambda,q}$. To establish the third property we only need to check that $\|x_n\|\to\|\pi(x)\|$ as $n\to+\infty$ for all $x\in\TT_{\lambda,q}^{(0)}$. Since the kernel of $\pi$ is $\K(\F_{\lambda,q})$, we have $\|\pi(x)\|=\limsup_{n\to+\infty}\|x_n\|$ for all $x\in\TT_{\lambda,q}^{(0)}$. On the other hand, by the previous lemma we have $\|\sigma_n(x_n)\|\to\|\pi(x)\|$. As the maps $\sigma_n$ are contractive, it follows that $\|\pi(x)\|\le \liminf_{n\to+\infty}\|x_n\|$. Hence $\|x_n\|\to\|\pi(x)\|$.
\ep

For $q=1$, the existence of a continuous field structure with fibers as in the theorem has been unconditionally shown by Landsman~\cite{MR1656992}. Namely, consider the $G$-equivariant ucp maps
$$
\breve\sigma_n\colon C(K^S\backslash G)\to \B(V_{n\lambda}),\quad\breve\sigma_n(f)=(\dim V_{n\lambda})\int_G f(g)U_n(g^{-1})P^h_{n\lambda}U_n(g)dg.
$$
Let us also write $\breve\sigma_\infty(f)$ for $f$. It is shown in~\cite{MR1656992}*{Theorem~1} that there is a unique continuous field of C$^*$-algebras with fibers as in Theorem~\ref{thm:cont-field} such that the sections $(\breve\sigma_n(f))_{n\in\Z_+\cup\{\infty\}}$ are continuous for all $f\in  C(K^S\backslash G)$.

\begin{proposition}
For $q=1$, the continuous field structure defined in Theorem~\ref{thm:cont-field} coincides with that defined in~\cite{MR1656992}. In other words, assuming that condition~\eqref{eq:conj1} is satisfied, the C$^*$-algebra $\TT^{(0)}_\lambda$ is generated by $\K(\F_\lambda)\cap \TT^{(0)}_\lambda$ and the elements $(\breve\sigma_n(f))_{n\ge0}$ for $f\in  C(K^S\backslash G)$.
\end{proposition}

\bp
It suffices to show that for all $x\in\TT^{(0)}_\lambda$ we have $\|x_n-\breve\sigma_n(\pi(x))\|\to0$ as $n\to+\infty$. Since the maps~$\breve\sigma_n$ form an asymptotic homomorphism by~\cite{MR1656992}*{Theorem~2}, it is enough to check this on elements of $\TT^{(0)}_\lambda$ that generate $\TT^{(0)}_\lambda$ modulo the compacts. As such elements we can take $S_\xi S_\zeta^*$ for $\xi,\zeta\in V_\lambda$, since $s_\xi s_\eta^*$ generate~$\OO_\lambda^{(0)}$ as a C$^*$-algebra, which follows from the facts that~$\OO_\lambda$ is abelian and the monomials in $s_\xi$ and $s_\eta^*$ containing the same number of generators and their adjoints span a dense subspace of $\OO_\lambda^{(0)}$.

We will use the same strategy as in the proof of Proposition~\ref{prop:star-commute}. For each $n$, consider two $G$-equivariant linear maps $A_n,B_n\colon V_{\lambda}\otimes\bar V_\lambda\to \B(V_{n\lambda})$ defined by
$$
A_n(\xi\otimes\bar\zeta)=S_\xi S_\zeta^*,\qquad B_n(\xi\otimes\bar\zeta)=\breve\sigma_n(\pi(S_\xi S_\zeta^*)).
$$
We need to show that $\|A_n-B_n\|\to0$. As in the proof of Proposition~\ref{prop:star-commute}, condition~\eqref{eq:conj1} implies that it suffices to check that  $\|(A_n(\xi\otimes\bar\zeta)-B_n(\xi\otimes\bar\zeta))\xi_{n\lambda}\|\to0$.
We know how $A_n(\xi\otimes\bar\zeta)\xi_{n\lambda}$  approximately looks like, so we need to show that $B_n(\xi_\lambda\otimes\bar\xi_\lambda)\xi_{n\lambda}$ becomes close to $\xi_{n\lambda}$ as $n\to+\infty$, while if $\xi,\zeta\in V_\lambda$ are weight vectors of norm one and either $\wt(\xi)\ne\lambda$ or $\wt(\zeta)\ne\lambda$, then $B_n(\xi\otimes\bar\zeta)\xi_{n\lambda}$ becomes close to zero.

Let us first consider $(B_n(\xi\otimes\bar\zeta)\xi_{n\lambda},\xi_{n\lambda})=(\breve\sigma_n(\pi(S_\xi S_\zeta^*))\xi_{n\lambda},\xi_{n\lambda})$. A simple computation shows that
$$
(\breve\sigma_n(f)\xi_{n\lambda},\xi_{n\lambda})=(\dim V_{n\lambda})\int_G f(g)|(U_n(g)\xi_{n\lambda},\xi_{n\lambda})|^2dg.
$$
This expression is known to converge to $f(e)$ as $n\to+\infty$ for all $f\in C(K^S\backslash G)$, see, e.g.,~\cite{MR2055928}*{Lemma~3.3} for a short proof. Recalling that $\pi(S_\xi)=C^\lambda_{\xi_\lambda,\xi}$, we get that
$$
(\breve\sigma_n(\pi(S_\xi S_\zeta^*))\xi_{n\lambda},\xi_{n\lambda})\to (\xi,\xi_\lambda)\overline{(\zeta,\xi_\lambda)}.
$$

In particular, $(\breve\sigma_n(\pi(S_{\xi_\lambda} S_{\xi_\lambda}^*))\xi_{n\lambda},\xi_{n\lambda})\to1$. As $\breve\sigma_n(\pi(S_{\xi_\lambda} S_{\xi_\lambda}^*))$ is a contraction, it follows that $\|\breve\sigma_n(\pi(S_{\xi_\lambda} S_{\xi_\lambda}^*))\xi_{n\lambda}-\xi_{n\lambda}\|\to0$. Similarly, if $\xi\in V_\lambda$ is a weight vector with $\wt(\xi)\ne\lambda$, we get $(\breve\sigma_n(\pi(S_{\xi} S_{\xi}^*))\xi_{n\lambda},\xi_{n\lambda})\to0$ and then conclude that $\|\breve\sigma_n(\pi(S_\xi S_\xi^*))\xi_{n\lambda}\|\to0$. Finally, if $\xi,\zeta\in V_\lambda$ are weight vectors of norm one and $\wt(\xi)\ne\lambda$, then
$$
\breve\sigma_n(\pi(S_\zeta S_\xi^*))^*\breve\sigma_n(\pi(S_\zeta S_\xi^*))\le\breve\sigma_n(\pi(S_\xi S_\zeta^*S_\zeta S_\xi^*))\le\breve\sigma_n(\pi(S_\xi S_\xi^*))
$$
and, since $\pi(S_\xi S_\zeta^*)=\pi(S_\zeta^*S_\xi)$,
$$
\breve\sigma_n(\pi(S_\xi S_\zeta^*))^*\breve\sigma_n(\pi(S_\xi S_\zeta^*))\le \breve\sigma_n(\pi(S_\xi^* S_\xi))=\breve\sigma_n(\pi(S_\xi S_\xi^*)),
$$
hence both $\|\breve\sigma_n(\pi(S_\zeta S_\xi^*))\xi_{n\lambda}\|$ and $\|\breve\sigma_n(\pi(S_\xi S_\zeta^*))\xi_{n\lambda}\|$ converge to zero as $n\to+\infty$. This completes the proof of the proposition.
\ep

\subsection{Compactification of discrete quantum spaces}

We will now consider the case $q>1$ and give another description of $\TT^{(0)}_{\lambda,q}$. Throughout this subsection we will assume that there exist $C>0$ and $0<t<1$ such that condition~\eqref{eq:conjS} is satisfied. Once again recall that this is satisfied for $G=SU(N)$, $\lambda\in P_{++}$ and $\lambda\in\mathbb N\omega_1$.

Define a discrete quantum set $\Gamma_{\lambda,q}$ by letting
$$
\ell^\infty(\Gamma_{\lambda,q})=\ell^\infty\text{-}\bigoplus^\infty_{n=0}\B(V_{n\lambda}).
$$
In other words, $\ell^\infty(\Gamma_{\lambda,q})$ consists of norm-bounded sequences of operators $x_n\in \B(V_{n\lambda})$. We view this algebra as a subalgebra of $\B(\F_{\lambda,q})$. We also define $c_0(\Gamma_{\lambda,q})\subset \ell^\infty(\Gamma_{\lambda,q})$ as the subalgebra of sequences $(x_n)_n$ such that $\|x_n\|\to0$. Equivalently, $c_0(\Gamma_{\lambda,q})=\ell^\infty(\Gamma_{\lambda,q})\cap\K(\F_{\lambda,q})$. Any unital C$^*$-subalgebra $C(\bar\Gamma_{\lambda,q})\subset \ell^\infty(\Gamma_{\lambda,q})$ containing $c_0(\Gamma_{\lambda,q})$ can be viewed as the algebra of continuous functions on a compactification of $\Gamma_{\lambda,q}$ and then the quotient $C(\partial\Gamma_{\lambda,q})=C(\bar\Gamma_{\lambda,q})/c_0(\Gamma_{\lambda,q})$ describes the corresponding boundary.

\smallskip

We are going to construct a compactification of $\Gamma_{\lambda,q}$ generalizing the construction of Vaes and Vergnioux for $SU_q(2)$ and other free orthogonal quantum groups~\cite{MR2355067}. It can be viewed as a quantum analogue of the end compactification of a tree.

Recall that $f_n$ denotes the projection $V_\lambda^{\otimes n}\to V_{n\lambda}$. Consider the inductive system of $G_q$-equivariant ucp maps
$$
\psi_{n,n+k}\colon \B(V_{n\lambda})\to \B(V_{(n+k)\lambda}),\quad T\mapsto f_{n+k}(T\otimes1)f_{n+k},
$$
where we use that $V_{(n+k)\lambda}\subset V_{n\lambda}\otimes V_{k\lambda}$, and define
\begin{align}
  C(\bar\Gamma_{\lambda,q}) & =\overline{\{x\in \ell^\infty(\Gamma_{\lambda,q})\mid \psi_{n,n+k}(x_n)=x_{n+k}\ \text{for all}\ n\ \text{large enough and}\ k\ge0\}}\notag\\
  & =\{x\in \ell^\infty(\Gamma_{\lambda,q})\mid \lim_{n\to+\infty}\sup_{k\ge0}\|\psi_{n,n+k}(x_n)-x_{n+k}\|=0\}.\label{eq:compactification}
\end{align}
It is clear that $C(\bar\Gamma_{\lambda,q})$ is an operator system containing $c_0(\Gamma_{\lambda,q})$, but it is not at all obvious that it is closed under multiplication.


\begin{theorem}\label{thm:gauge-invariant-sub}
Let $G$ be a simply connected semisimple compact Lie group, $\lambda\in P_+\setminus\{0\}$ and $q>1$. Let $S=\{\alpha\in\Pi: (\lambda,\alpha)=0\}$. Assume condition~\eqref{eq:conjS} is satisfied for some $C>0$ and $0<t<1$. Then the operator system $ C(\bar\Gamma_{\lambda,q})$ defined by~\eqref{eq:compactification} is a C$^*$-subalgebra of $\ell^\infty(\Gamma_{\lambda,q})\subset \B(\F_{\lambda,q})$, and we have $C(\bar\Gamma_{\lambda,q})=\TT^{(0)}_{\lambda,q}$. Moreover, the surjective $G_q$-equivariant homomorphism $\TT_{\lambda,q}\to C( G^\lambda_q\backslash G_q)$ from Theorem~\ref{thm:main} gives rise to an isomorphism
$$
C(\partial\Gamma_{\lambda,q})=C(\bar\Gamma_{\lambda,q})/c_0(\Gamma_{\lambda,q})\cong C(K^S_q\backslash G_q).
$$
\end{theorem}

This generalizes \cite{MR4705666}*{Theorem 4.1} for $G_q=SU_q(2)$. Our strategy will be the same as in~\cite{MR4705666}, but several arguments require some adjustments.

\smallskip

The following estimate will play a key role in the proof.

\begin{lemma}[{cf.~\cite{MR4705666}*{Lemma~4.2}}]\label{f-estimate}
For all $n\ge2$, we have
$$
\|f_{n+1}-(f_n\otimes1)(1\otimes f_n)\|\le \|(1-1\otimes P^h_{n\lambda})P^h_{\lambda,n\lambda}\|+\big\|f_{\lambda,(n-1)\lambda}\big((1-P^h_\lambda)\otimes P^h_{(n-1)\lambda}\big)\big\|.
$$
\end{lemma}

\bp
Since $f_{n+1}$ is dominated by $1\otimes f_n$ and $f_n\otimes1$, the inequality in the statement of the lemma is equivalent~to
$$
\|(f_n\otimes1)|_{(V_\lambda\otimes V_{n\lambda})\ominus V_{(n+1)\lambda}}\|\le \|(1-1\otimes P^h_{n\lambda})P^h_{\lambda,n\lambda}\|+\big\|f_{\lambda,(n-1)\lambda}\big((1-P^h_\lambda)\otimes P^h_{(n-1)\lambda}\big)\big\|.
$$

As we already used in the proof of Proposition~\ref{prop:star-commute}, since $f_n\otimes 1$ is a morphism, its norm is determined by the restriction to the subspace of vectors killed by the~$E_i$'s. The projection onto the subspace of such vectors in $(V_\lambda\otimes V_{n\lambda})\ominus V_{(n+1)\lambda}$ is $P^h_{\lambda,n\lambda}-(P^h_\lambda)^{\otimes(n+1)}$, where we view $P_{\lambda,n\lambda}$ as an operator on~$V^{\otimes(n+1)}_\lambda$. Therefore
\begin{multline}\label{eq:fn}
\|(f_n\otimes1)|_{(V_\lambda\otimes V_{n\lambda})\ominus V_{(n+1)\lambda}}\| = \|(f_n\otimes1)(P^h_{\lambda,n\lambda}-(P^h_\lambda)^{\otimes(n+1)})\|\\
\le \|(f_n\otimes1)\big((1\otimes P^h_{n\lambda})P^h_{\lambda,n\lambda}-(P^h_\lambda)^{\otimes(n+1)}\big)\|+\|(1-1\otimes P^h_{n\lambda})P^h_{\lambda,n\lambda}\|.
\end{multline}
We can write
$$
(1\otimes P^h_{n\lambda})P^h_{\lambda,n\lambda}-(P^h_\lambda)^{\otimes(n+1)}=\big((1-P^h_\lambda)\otimes P^h_{(n-1)\lambda}\otimes P^h_\lambda\big)P^h_{\lambda,n\lambda}.
$$
Since $f_n$ is nothing other than $f_{\lambda,(n-1)\lambda}$ viewed as an operator on $V_\lambda^{\otimes n}$, it follows that
$$
\|(f_n\otimes1)\big((1\otimes P^h_{n\lambda})P^h_{\lambda,n\lambda}-(P^h_\lambda)^{\otimes(n+1)}\big)\|\le \big\|f_{\lambda,(n-1)\lambda}\big((1-P^h_\lambda)\otimes P^h_{(n-1)\lambda}\big)\big\|.
$$
This, taken together with~\eqref{eq:fn}, proves the lemma.
\ep

Using this lemma and the assumption~\eqref{eq:conjS} the next several results are proved identically to~\cite{MR4705666}, so we omit the proofs.
But first we need to introduce some notation. Similarly to the left creation operators $S_\xi$ we can define right creation operators $R_\xi$ on $\F_{\lambda,q}$:
$$
R_\xi\eta=f_{n+1}(\eta\otimes\xi)\quad\text{for}\quad \xi\in V_\lambda\ \  \text{and}\ \ \eta \in V_{n\lambda}.
$$

\begin{lemma}[{\cite{MR4705666}*{Lemma~4.3}}]\label{lem:commutators}
For all $\xi,\zeta\in V_\lambda$, we have $[S_\xi,R_\zeta]=0$. There is a constant $c>0$ such that
$$
\|[S_\xi^*,R_\zeta]|_{V_{n\lambda}}\|\le ct^n\|\xi\|\,\|\zeta\|\quad\text{for all}\quad\xi,\zeta\in V_\lambda\quad\text{and}\quad n\ge0.
$$
\end{lemma}

This immediately gives the following result.

\begin{corollary}\label{cor:commutation-mod-compacts}
For every $S\in\TT_{\lambda,q}$ and every $R$ in the C$^*$-algebra generated by the operators~$R_\xi$, $\xi\in V_\lambda$, we have
$[S,R]\in\K(\F_{\lambda,q})$.
\end{corollary}

Fix an orthonormal basis $(\xi_i)_i$ in $V_\lambda$ and write $R_i$ for $R_{\xi_i}$. Similarly to~\eqref{eq:e-0}, we have $\sum^N_{i=1}R_iR_i^*=1-e_0$. Define a contractive cp map
$$
\Theta\colon \B(\F_{\lambda,q})\to \B(\F_{\lambda,q}),\quad \Theta(T)=\sum_iR_iTR_i^*.
$$
Since $\K(\F_{\lambda,q})\subset\TT_{\lambda,q}$, by Corollary~\ref{cor:commutation-mod-compacts} this map leaves $\TT_{\lambda,q}$ globally invariant. It also leaves $\ell^\infty(\Gamma_{\lambda,q})$ globally invariant and we have
\begin{equation}\label{eq:theta}
\Theta^k(x)_{n+k}=\psi_{n,n+k}(x_n)\quad\text{for all}\quad x\in \ell^\infty(\Gamma_{\lambda,q}),\quad n,k\ge0.
\end{equation}

Denote by $\A_{\lambda,q}$ the unital $*$-subalgebra of $\TT_{\lambda,q}$ generated by the elements $S_\xi$, $\xi\in V_\lambda$. Let $\A^{(0)}_{\lambda,q}$ be the gauge-invariant part of $\A_{\lambda,q}$. It is spanned by the monomials in $S_\xi$ and $S_\xi^*$ containing an equal number of creation and annihilation operators. Using Lemma~\ref{lem:commutators} one gets the following result.

\begin{lemma}[{\cite{MR4705666}*{Lemma~4.6}}]\label{lem:theta-estimate}
For every $x\in\A^{(0)}_{\lambda,q}$, there exists a constant $C_x>0$ such that
$$
\|x_{n+k}-\Theta^k(x)_{n+k}\|\le C_xt^n \quad\text{for all}\quad n,k\ge0.
$$
\end{lemma}

\bp[Proof of Theorem~\ref{thm:gauge-invariant-sub}]
The last part of the theorem follows from Lemma~\ref{lem:flag}, so we only need to establish the equality $C(\bar\Gamma_{\lambda,q})=\TT^{(0)}_{\lambda,q}$.

From~\eqref{eq:theta} and Lemma~\ref{lem:theta-estimate} we see that $\A^{(0)}_{\lambda,q}\subset C(\bar \Gamma_{\lambda,q})$. Hence $\TT^{(0)}_{\lambda,q}\subset C(\bar \Gamma_{\lambda,q})$.
In order to prove the opposite inclusion, take any $x\in\ell^\infty(\Gamma_{\lambda,q})$ such that $$x_{n_0+k}=\psi_{n_0,n_0+k}(x_{n_0})$$ for some $n_0$ and all $k\ge0$. Since the maps $\psi_{n_0,n_0+k}$ are $G_q$-equivariant with respect to the adjoint action, we may assume that $x_{n_0}$ lies in an isotypical component $\B(V_{n_0\lambda})_\nu\subset \B(V_{n_0\lambda})$ with highest weight $\nu\in P_+$.

Denote also by $C(K^S_q\backslash G_q)_\nu\subset\C[K^S_q\backslash G_q]$ the isotypical component (with respect to the action by right translations) with highest weight $\nu$. Its multiplicity is equal to the dimension of the space of $K^S_q$-invariant vectors in $V_\nu$. By Proposition~\ref{prop:mult} this dimension is also not smaller than the multiplicity of the isotypical component $\B(V_{n\lambda})_\nu$, since this multiplicity equals the multiplicity of $V_{n\lambda}$ in $V_\nu\otimes V_{n\lambda}$. Thus,
\begin{equation}\label{eq:iso-dim}
\dim \B(V_{n\lambda})_\nu\le \dim C(K^S_q\backslash G_q)_\nu  \quad\text{for all}\quad n\ge0.
\end{equation}

Choose a lift of the embedding map $C(K^S_q\backslash G_q)_\nu\to C(K^S_q\backslash G_q)$ to a $G_q$-equivariant linear map
$$
\rho\colon  C(K^S_q\backslash G_q)_\nu\to \A^{(0)}_{\lambda,q},
$$
so that $\pi(\rho(a))=a$ for $a\in  C(K^S_q\backslash G_q)_\nu$. By Lemma~\ref{lem:theta-estimate} applied to the operators in the image of $\rho$, we can find $c>0$ such that
\begin{equation}\label{eq:pi-theta}
\|\rho(a)_{n+k}-\Theta^k(\rho(a))_{n+k}\|\le ct^n\|a\|\quad\text{for all}\quad a\in C(K^S_q\backslash G_q)_\nu\quad\text{and}\quad n,k\ge0.
\end{equation}

Next, take any $0<\delta<1$. By Lemma~\ref{lem:norm-convergence} we can find $n_1\ge0$ such that $\|\rho(a)_n\|\ge \delta\|a\|$ for all $a\in C(K^S_q\backslash G_q)_\nu$ and $n\ge n_1$. It follows that for all $n\ge n_1$ the map
$$
C(K^S_q\backslash G_q)_\nu\to \B(V_{n\lambda})_\nu,\quad a\mapsto\rho(a)_n,
$$
is a linear isomorphism.

Now, fix $\eps>0$ and choose $n\ge\max\{n_0,n_1\}$ such that $ct^n\|x\|<\eps \delta$. For the unique $a\in C(K^S_q\backslash G_q)_\nu$ such that $\rho(a)_n=x_n$ we have $\|a\|\le \delta^{-1}\|\rho(a)_n\|\le \delta^{-1}\|x\|$. For all $k\ge0$ we have
$$
x_{n+k}=\Theta^k(x)_{n+k}=\Theta^k(\rho(a))_{n+k}.
$$
Hence, applying again~\eqref{eq:pi-theta}, we get
$$
\|x_{n+k}-\rho(a)_{n+k}\|\le ct^n\|a\|<\eps.
$$
Therefore, modulo the compacts, $x$ is $\eps$-close to $\A^{(0)}_{\lambda,q}$. Hence $x\in\TT^{(0)}_{\lambda,q}$.
\ep

\bigskip

\appendix
\section{Asymptotics of Clebsch--Gordan coefficients for \texorpdfstring{$SU_q(N)$}{SUq(N)}}\label{appendix}

The goal of this appendix is to sketch a proof of Theorem~\ref{thm:SUN}. As discussed in Section~\ref{ssec:conj}, it suffices to consider $q\ge1$.

\subsection{Regular weights}
We first consider the regular weights. We want to apply Lemma~\ref{lem:regular} to $V=V_{\omega_1}$. Thus, we need to verify that
\begin{equation}\label{eq:app-conj}
\|P^h_{\omega_1,\mu}-1\otimes P^h_\mu\|\to0\quad\text{as}\quad N(\mu)\to+\infty,
\end{equation}
and the convergence is exponentially fast when $q>1$.

Recall that every highest weight $\mu$ is described by a tuple $(\mu_1,\dots,\mu_{N-1},0)$, where the $\mu_i$'s are integers satisfying $\mu_1 \geq \cdots \geq \mu_{N-1} \geq 0$. Then $\mu(i)=\mu_i-\mu_{i+1}$ (with $\mu_N=0$) and therefore
$$
N(\mu)=\min_{1\le i\le N-1}(\mu_i-\mu_{i+1}).
$$

The Clebsch--Gordan coefficients for $V_{\omega_1}\otimes V_\mu$ are described in terms of the Gelfand--Tsetlin bases. Namely, the basis elements of $V_\mu$ are parameterized by the following arrays:
\[
\mathbf{r}=
\begin{pmatrix}
r_{11} & \dots & r_{1,N-1} & 0\\
r_{21} & \dots & r_{2,N-1} &\\
\vdots & & & \\
r_{N1} & & &
\end{pmatrix},
\]
where the $r_{ij}$'s are integers satisfying $r_{1j}=\mu_j$ for $j=1,\dots,N-1$ and $r_{ij} \geq r_{i+1,j} \geq r_{i,j+1} \geq 0$ for all $i,j$. In this basis the highest weight vector $\xi_\mu$ can be taken to be
\[
\mathbf{r}(\mu)=
\begin{pmatrix}
\mu_{1} & \dots & \mu_{N-1} & 0\\
\mu_{1} & \dots & \mu_{N-1} &\\
\vdots & & & \\
\mu_{1} & & &
\end{pmatrix}.
\]
The elements $e_i$ of the standard basis of $V_{\omega_1}=\C^N$ correspond to the arrays $\mathbf{e}^i=(e^i_{ab})$ compactly written as
\[
e^i_{ab}=
\begin{cases}
1, &\text{if}\ \ 1 \leq a \leq N-i+1\ \ \text{and}\ \ b=1,\\
0, &\text{otherwise}.
\end{cases}
\]

Now, for $N(\mu)\ge1$ we have
\[
V_{\omega_1} \otimes V_{\mu} \cong \bigoplus_{i=1}^N V_{\mu^i},
\]
where $\mu^i=(\mu_1,\dots,\mu_i+1,\dots,\mu_{N-1},0)$ for $i=1,\dots,N-1$ and $\mu^N=(\mu_1-1,\dots,\mu_{N-1}-1,0)$. In order to prove~\eqref{eq:app-conj} it then suffices to show that for all $1\le i\le N$ we have
$$
|(\mathbf{e}^i \otimes \mathbf{r}(\mu),\mathbf{r}(\mu^i))|\to 1\quad\text{as}\quad N(\mu)\to+\infty,
$$
and the convergence is exponentially fast when $q>1$. This is indeed true, since up to phase factors we have
\begin{equation}\label{eq:CG}
(\mathbf{e}^i \otimes \mathbf{r}(\mu),\mathbf{r}(\mu^i)) = q^{\frac{1}{2}(i-1)} \sqrt{\frac{\prod_{j=1}^{i-1}[\mu_{j}-\mu_{i}-j+i-1]_q}{\prod_{j=1}^{i-1}[\mu_{j}-\mu_{i}-j+i]_q}},
\end{equation}
which can be deduced from formulas for the Clebsch--Gordan coefficients of $V_{\omega_1}\otimes V_\mu$, see, e.g., \cite{klimyk-schmudgen:book-quantum-groups}*{p. 220} and \cite{chakraborty-pal:spectral-triples}*{Section 4.2}. Note that the conventions in these references are different from ours, so $q$ in formulas there should be replaced by $q^{-1}$.

\subsection{Multiples of \texorpdfstring{$\omega_1$}{omega1}}
Next we consider weights of the form $\lambda=m\omega_1$. The first convergence in~\eqref{eq:conj1} follows again from known formulas for the Clebsch--Gordan coefficients, see, e.g.,~\cite{alisaukas-smirnov:clebsch-gordan}, but their general form is significantly more complicated for $m\ge2$ compared to the formulas used above, so it is easier to argue as follows.

For $N=2$ the convergence holds by the case of regular weights. Assume therefore that $N\ge3$. We identify $V_{n\omega_1}$ with the space $\C_q[e_1,\dots,e_N]_n$ of homogeneous polynomials of degree~$n$ in $\C_q[e_1,\dots,e_N]$. For $n\ge m$, we have
$$
V_{m\omega_1}\otimes V_{n\omega_1}\cong\bigoplus^m_{k=0}V_{\mu^{n,k}},
$$
where $\mu^{n,k}=(n+k,m-k,0,\dots,0)$. Take any $0\le k\le m$. Let $\xi^{n,k}\in\C_q[e_1,e_2]_m\otimes \C_q[e_1,e_2]_n\subset \C_q[e_1,\dots,e_N]_m\otimes \C_q[e_1,\dots,e_N]_n$ be the unique up to a phase factor unit vector of weight $(n+k,m-k,0,\dots,0)$ killed by~$E_1$, which exists by the fusion rules for $SU_q(2)$. It is obvious that it is also killed by $E_i$ for $2\le i\le N-1$. Therefore this is a highest weight vector determining the component $V_{\mu^{n,k}}$ of the tensor product $V_{m\omega_1}\otimes V_{n\omega_1}$. Consider also the unique up to a phase factor unit vector $\zeta^k\in\C_q[e_1,e_2]_m\subset\C_q[e_1,\dots,e_N]_m$ of weight $(k,m-k,0,\dots,0)$. Explicitly, we can take
$$
\zeta^{k}=q^{-k(m-k)/{2}}\Big(\frac{[m]_q!}{[k]_q![m-k]_q!}\Big)^{1/2}e_1^ke_2^{m-k},
$$
where we used~\eqref{eq:inner-products} to get the correct normalization. By the case $N=2$ we know that
$$
|(\zeta^k\otimes e_1^n,\xi^{n,k})|\to 1\quad\text{as}\quad n\to+\infty,
$$
and the convergence is exponentially fast when $q>1$. Therefore the first convergence in~\eqref{eq:conj1} holds for $\lambda=m\omega_1$.

For the second convergence in~\eqref{eq:conj1}, consider $\lambda=\omega_1$. For $2\le i\le N$, using again~\eqref{eq:inner-products} we~get
$$
\|f_{\omega_1,n\omega_1}(e_i\otimes\xi_{n\omega_1})\|=\|e_ie_1^n\|=q^{-n}\|e_1^ne_i\|=q^{-{n}/{2}}\frac{1}{\sqrt{[n+1]_q}}.
$$
Therefore we have $\big\|f_{\omega_1,n\omega_1}\big((1-P^h_{\omega_1})\otimes P^h_{n\omega_1}\big)\big\|\to0$ as $n\to+\infty$, and the convergence is exponentially fast when $q>1$. A similar argument works also for $\lambda=m\omega_1$.

\bigskip

\end{document}